\documentclass[11pt]{article}

\usepackage[utf8]{inputenc} 
\usepackage{booktabs}       
\usepackage{amsfonts}       
\usepackage{nicefrac}       
\usepackage{microtype}      
\usepackage{amsmath}
\usepackage[usenames, dvipsnames]{color}
\usepackage{appendix}
\usepackage{comment}
\usepackage{fullpage}
\usepackage{hyperref}
\usepackage{natbib}

\usepackage{authblk}
\usepackage{xr}

\newcommand{\JG}[1]{{\color{magenta}JG: #1}}

\newtheorem{theorem}{Theorem}
\newtheorem{lemma}[theorem]{Lemma}

\begin{document}

\title{Posterior Contraction Rates for Gaussian Cox Processes with Non-identically Distributed Data}

\author[1,3]{James A. Grant\thanks{j.grant@lancaster.ac.uk; corresponding author}}
\author[2,3]{David S. Leslie\thanks{d.leslie@lancaster.ac.uk}}
\affil[1]{STOR-i Centre for Doctoral Training, Lancaster University, UK}
\affil[2]{Department of Mathematics and Statistics, Lancaster University, UK}
\affil[3]{PROWLER.io, Cambridge, UK}
\maketitle

\begin{abstract}
This paper considers the posterior contraction of non-parametric Bayesian inference on non-homogeneous Poisson processes. We consider the quality of inference on a rate function $\lambda$, given non-identically distributed realisations, whose rates are transformations of $\lambda$. Such data arises frequently in practice due, for instance, to the challenges of making observations with limited resources or the effects of weather on detectability of events. We derive contraction rates for the posterior estimates arising from the  Sigmoidal Gaussian Cox Process and Quadratic Gaussian Cox Process  models. These are popular models where $\lambda$ is modelled as a logistic and quadratic transformation of a Gaussian Process respectively. Our work extends beyond existing analyses in several regards. Firstly, we consider non-identically distributed data, previously unstudied in the Poisson process setting. Secondly, we consider the Quadratic Gaussian Cox Process model, of which there was previously little theoretical understanding. Thirdly, we provide rates on the shrinkage of both the width of balls around the true $\lambda$ in which the posterior mass is concentrated and on the shrinkage of posterior mass outside these balls - usually only the former is explicitly given. Finally, our results hold for certain finite numbers of observations, rather than only asymptotically, and we relate particular choices of hyperparameter/prior to these results. 
\end{abstract}

\section{Introduction} \label{sec::intro}

The non-homogeneous Poisson process (NHPP) is the most widely used model for inference on point process data. It is parameterised by a non-negative rate function $\lambda$ and satisfies the key property that the expected number of events in any area is equal to the integral of the rate function over that area. A Gaussian Cox process (GCP) is a nonparametric Bayesian version of the NHPP model where $\lambda$ is modelled as a transformation of a Gaussian process (GP). In this paper we consider two classes of GCP, the Sigmoidal GCP (SGCP) of \citep{AdamsEtAl2009} and the Quadratic GCP (QGCP) of \citep{LloydEtAl2015}. In the SGCP the rate function is modelled as a multiple of a logistic transformation of a GP. In the QGCP the rate function is modelled as the square of a GP. This paper is concerned with the quality of posterior inference on $\lambda$ arising from these models. Specifically we are interested in the rate at which the expected posterior mass the models assign to functions far from the true $\lambda$ decreases. 

The GCP is a model over functions and is defined on some space of non-negative functions $\Lambda$. Given a true rate function $\lambda_0 \in \Lambda$, observed data $X_{1:n}$ collected over $n \in \mathbb{N}$ timesteps, and a relevant distance $d_n(\lambda,\lambda')$ defined for all $\lambda,\lambda' \in \Lambda$, we look for results of the form \begin{equation}
\mathbb{E}_{\lambda_0}\bigg(\Pi\big(\lambda \in \Lambda : d_n(\lambda,\lambda_0)\geq \epsilon_n | X_{1:n}\big) \bigg) \leq f_n \label{eq::example_contraction}
\end{equation} for decreasing sequences $\epsilon_n, f_n$, where $\Pi(\cdot|X_{1:n})$ denotes the posterior probability mass and $\mathbb{E}_{\lambda_0}$ denotes expectation with respect to the probability measure implied by $\lambda_0$. The sequences $\epsilon_n, f_n$ define the rate of posterior contraction of a model. If such a bound holds for certain $\epsilon_n, f_n \rightarrow 0$ as $n \rightarrow \infty$ this displays that the model is consistent. However, we are also interested in the order of the sequences and for which $n$ results of the form \eqref{eq::example_contraction} can be identified.

Asymptotic consistency results of the form \begin{equation*}
\mathbb{E}_{\lambda_0}\bigg(\Pi\big(\lambda \in \Lambda : d_n(\lambda,\lambda_0)\geq \epsilon_n | X_{1:n}\big) \bigg) \rightarrow 0
\end{equation*} as $n \rightarrow \infty$ are prevalent in the Bayesian nonparametrics  literature; for example \citep{KirichenkoVanZanten2015} gives such a result for i.i.d. $X_{1:n}$ under the SGCP and a broader family of GCPs which have a smooth and bounded link function. Such asymptotic results are undoubtedly useful contributions to the understanding Bayesian models and inference, however they provide limited support to finite-time analyses thereof. We extend beyond existing results in four important regards by \begin{enumerate}
\item Providing results for independent non-identically distributed (i.n.i.d.) data,
\item Providing results for the QGCP model as well as the SGCP,
\item Providing a rate on the shrinkage of the posterior mass $f_n$ (as well as on $\epsilon_n$), and
\item Providing results for finite values of $n$, not only asymptotically, and relating specific choices of hyperpriors (and parameters) to these results.
\end{enumerate}

Studying i.n.i.d. data is in contrast to the majority of previous studies of non-parametric inference on NHPPs. However, ours is an important, practically-relevant setting. Commonly when observing point process data, the detection of events may be imperfect. This may be due to visibility conditions, unreliable signals or the fallability of observation equipment. A result of this is that while events may occur independently and according to a stationary process, the distribution of \emph{observed} events can vary as data is collected. Equally, different subsections of a region of interest may be observed at different rates by design. Data collectors may be more readily able to gather data in a particular region, resources may be too costly to gather the same quality of information everywhere or multiple sub-investigations may be combined to form a joint dataset. As the GCP models are typically used to model situations with underlying spatial smoothness and covariance structure, a unified analysis is still desirable, however existing contraction results only handle the setting where an entire region of interest has been observed uniformly. The results we obtain in this paper apply to the setting where (whether through design or imprecision) different rates of observation have been applied at different locations. Therefore, we present results that are more relevant to the practical settings in which GCP models are utilised than those which consider only identically distributed data.

The QGCP model has recently received attention in the literature \citep{LloydEtAl2015, JohnHensman2018} as a model for NHPP inference, due to the ability to carry out fast and accurate inference. Previously however, there was little theoretical understanding of the model. We provide theoretical foundations for this new variant of the GCP model. This is non-trivial since the link function in the QGCP is not bounded, in contrast with the SGCP. Consequently, we find the rate of contraction to be lower for the QGCP than for the traditional SGCP. 

Providing a rate $f_n$ on the shrinkage of the posterior mass for finite values of $n$ is also an important development.  A trend in the existing literature is to focus on the asymptotic results and present that if the width of a ball around the true rate is chosen to decrease at the correct rate (with respect to the number of observations) then the probability of lying outside this ball tends to 0 as the number of observations goes to infinity. Such results are typically cleaner, and clearly demonstrate the consistency of a method, while the finite-time result can usually be extracted from the proofs provided for such results if desired. If the rate $f_n$ is explicitly given or can be inferred, it is often only specified as holding for ``sufficiently large" $n$. Inferring the rate $f_n$ and determining the order of $n$ that qualifies as sufficiently large, can be challenging to users of these results. By explicitly giving a form of $f_n$ and quantifying the values of $n$ (in terms of functions of the chosen hyperparameters) for which it is valid, we present a more informative set of results that are useful for end-users of this theory.

One use case of these results is in the theory of sequential decision making problems. In sequential sensor placement problems such as that studied in \citep{GrantEtAl2019SP}, decision makers adaptively select intervals over which to observe events subject to costs on the length of the interval, with an overall aim of maximising a cumulative reward. To do so optimally, balancing between exploring undersampled regions and ``exploiting" - making repeated samples of areas where the rate is known to be high - is required. To understand the optimal balance of exploration and exploitation, one must understand the rate at which the inference model used contracts. Previous work on these problems has relied on assuming simpler inference models to obtain performance guarantees \citep{GrantEtAl2018UCB,GrantEtAl2019SP}. Guarantees on the contraction of Cox process posteriors \emph{with rates on the posterior mass} will be important in the design and analysis of more sophisticated approaches to these problems.

Another use for these results is in experimental design and resource planning problems. It is valuable for decision-makers to know the expected level of uncertainty in a rate function given a certain number of observations. They can then appropriately design sampling strategies or deploy resources to collect information in a way that is tailored to achieving a certain level of confidence in the inference.

In the remainder of this section, we discuss related work in GCPs and general contraction results for Bayesian models. In Section 2 we formally introduce our GCP models and notation. Section 3 includes all our main theoretical results and proofs, and in Section 4 we conclude with a discussion. Throughout we have aspired to make our assumptions transparent and demonstrate how they can be met. In Appendix H we verify that all assumed conditions can be satisfied for finite numbers of observations.

\subsection{Related Literature}

The Cox process \citep{Cox1955} is a class of doubly stochastic process where the rate function of an non-homogeneous Poisson process (see e.g. \citep{MollerWaagepetersen2003}) is modelled as another stochastic process. The Gaussian Cox process (GCP), as mentioned above, is a particular subset of this class where the rate function of the NHPP is modelled via a transformation of a Gaussian process \citep{RasmussenWilliams2006}. Three main transformations have been proposed yielding three main models. Firstly, the Log-Gaussian Cox Process (LGCP) of \citep{RathbunCressie1994} and \citep{MollerEtAl1998} where $\lambda$ is modelled as an exponential transformation of a GP. Secondly, the Sigmoidal-Gaussian Cox Process (SGCP) of \citep{AdamsEtAl2009} where $\lambda$ is modelled as a multiple of a logistic transformation of a GP. Finally, the Quadratic-Gaussian Cox Process (QGCP) of \citep{LloydEtAl2015} where $\lambda$ is modelled as a quadratic transformation of a GP. We focus on the SGCP and QGCP models, as (for reasons discussed fully in Section 4) the LGCP model requires separate techniques to derive a contraction result.

General results for the contraction of posterior density estimates given i.i.d. data are available thanks to the seminal papers \citep{GhosalEtAl2000} and \citep{GhosalVanDerVaart2001}. The link between density estimation and function estimation is exploited in \citep{BelitserEtAl2013} to extend this work to show contraction rates for Bayesian Poisson process inference subject to appropriate prior conditions. Furthermore, \citep{BelitserEtAl2013} proposes a spline based prior satisfying these conditions. The result of \citep{BelitserEtAl2013} and GP concentration results of \citep{VanDerVaartVanZanten2009}  are used by \citep{KirichenkoVanZanten2015} to show an asymptotic rate of posterior contraction for the SGCP - \citep{KirichenkoVanZanten2015} is the existing work most similar to our contribution. However we are able to move beyond i.i.d. data to the independent non-identically distributed (i.n.i.d.) case, thanks to the work of \citep{GhosalVanDerVaart2007} in deriving contraction results for posterior density estimates under such data.

\section{Model} \label{sec::model}

In this section we introduce the data generating model and two prior models considered in the paper, along with other relevant notation required to understand our main results.

\subsection{Likelihood}

We consider an NHPP with bounded non-negative rate function $\lambda_0$ on $[0,1]^d$. We suppose that $n$ independent realisations of the NHPP $\tilde{X}_1,...,\tilde{X}_n$ are generated. Each realisation $j$ consists of a collection $m_j$ of points $\{\tilde{X}_j^1,...,\tilde{X}_j^{m_j}\} \in [0,1]^d$. We write \begin{equation*}
\tilde{X}_j=\sum_{i=1}^{m_j} \delta_{\tilde{X}_j^i}, \quad j=1,...,n
\end{equation*} where $\delta_{x}$ denotes the Dirac measure at $x$. By the definition of the NHPP model, each realisation $j$ is distributed such that the number of points in any set $R \subseteq S$, denoted $\tilde{X}_j(R)$ follows a Poisson distribution with mean $\int_B \lambda(s)ds$. Furthermore $\tilde{X}_j(R_1), \tilde{X}_j(R_2)$ are independent if the sets $R_1,R_2 \subseteq S$ are disjoint. 

Under our model, the realisations $\tilde{X}_1,...,\tilde{X}_n$ are not directly observed. Instead, so-called \emph{filtered realisations} ${X}_{1:n}={X}_1,...,{X}_n$ are observed. The events in a filtered realisation ${X}_j$ are a subset of the events in the corresponding \emph{raw realisation} $\tilde{X}_j$. The relationship between $X_{1:n}$ and $\tilde{X}_{1:n}$ is governed by a set of \emph{filtering functions} $\gamma_{1:n}=\gamma_1,...,\gamma_n$. 

Each filtering function $\gamma_j: [0,1]^d \rightarrow [0,1]$ evaluated at a point $s \in S$ gives the probability of observing an event in $\tilde{X}_j$ given that it has occurred at location $s$. Every event that occurs in $\tilde{X}_j$ is observed or not independently according to these probabilities. 

By standard results, ${X}_j$ is distributed according to an NHPP with rate $\gamma_i\lambda_0$. That is to say, the $n$ filtered realisations ${X}_{i:n}$ are then realisations of  \emph{independent, non-identically distributed} NHPPs with rates $\gamma_1\lambda_0,...,\gamma_n\lambda_0$ respectively. 

It follows that the likelihood of a particular set of observations ${X}_{1:n}={X}_1,...,{X}_n$ given a rate function $\lambda$ and filtering functions $\gamma_{1:n}$ can be written \begin{displaymath}
\mathcal{L}({X}_{1:n}|\lambda_0,\gamma_{1:n}) = \prod_{j=1}^n \exp\bigg({\int_S  \gamma_j(s)\lambda_0(s)d{X}_j(s) - \int_S (\gamma_j(s)\lambda_0(s)-1) ds}\bigg),
\end{displaymath} using the law of the realisation ${X}_j$ as given by Proposition 6.1 of \citep{Karr1986}. We note that the case of i.i.d. data as considered in \citep{KirichenkoVanZanten2015} and \citep{GugushviliEtAl2018} is a special case of this model, where $\gamma_{j}(x)=1, \forall x \in [0,1]^d, \forall j=1,\dots,n$.

\subsection{Prior Models}

In this paper we consider two Bayesian models of the Poisson process where the rate function $\lambda_0$ is  modelled a priori as a transformation of a Gaussian process. Under the SGCP model \citep{AdamsEtAl2009}, the true rate function is modelled a priori as \begin{equation}
    \lambda(s) = \lambda^* \sigma(g(s))=\lambda^*(1+e^{-g(s)})^{-1} \quad s \in S
\end{equation}
where $\lambda^*>0$ is a scalar hyperparameter endowed with an independent Gamma prior and $g$ is a zero-mean GP.  The sigmoidal transformation $\sigma$ is bounded in $[0,1]$ so the hyperparameter $\lambda^*$ models the maximum of the rate function, $||\lambda_0||_\infty$. The QGCP model \citep{LloydEtAl2015} uses a more straightforward transformation. The rate function is modelled a priori as \begin{equation}
    \lambda(s) = (g(s))^2 \quad s \in S
\end{equation}
where again, $g$ is a GP.

For both models, we specify certain additional properties of the GP to support our subsequent analyses. These conditions are standard in the posterior contraction literature  \citep{VanDerVaartVanZanten2008,VanDerVaartVanZanten2009,KirichenkoVanZanten2015}. We require that the covariance kernel $f$ of the GP $g$, can be given in its spectral form by \begin{equation}
Ef(s)f(s') = \int e^{-i<\xi,l(s'-s)>}\mu(\xi)d\xi, \quad s,s' \in S. \label{eq::covariance}
\end{equation}
Here $l>0$ is an (inverse) length scale parameter and $\mu$ is a spectral density on $\mathbb{R}^d$ such that the map $a \mapsto \mu(a\xi)$ on $(0,\infty)$ is decreasing for every $\xi \in \mathbb{R}^d$ and that satisfies \begin{displaymath}
\int e^{\delta||\xi||}\mu(d\xi) < \infty
\end{displaymath}
for some $\delta >0$. Condition \eqref{eq::covariance} is satisfied, for instance, by the squared exponential covariance function \begin{equation*}
E f(s)f(s')= e^{-l^2||s-s'||^2}, \quad s, s' \in S
\end{equation*} since it corresponds to a centred Gaussian spectral density.

The length scale parameter should have a prior $\pi_l$ on $[0,\infty)$ which satisfies \begin{equation}
C_1x^{q_1}\exp(-D_1x^d\log^{q_2} x) \leq \pi_l(x) \leq C_2x^{q_1}\exp(-D_2x^d\log^{q_2} x) \label{eq::lengthscalecondition}
\end{equation}
for positive constants $C_1, C_2, D_1, D_2$, non-negative constants $q_1,q_2$, and every sufficiently large $x>0$. In particular if $l^d$ is endowed with a Gamma$(a,b)$ prior, then \begin{equation*}
    \pi_l(x)=\frac{b^a d}{\Gamma(a)}x^{da-1}\exp(-bx^d)
\end{equation*} for $x>0$, and thus \eqref{eq::lengthscalecondition} is satisfied with $C_1=C_2=\frac{b^a d}{\Gamma(a)}$, $D_1=D_2=b$, $q_1=da-1$, and $q_2=0$. We will assume a Gamma prior on $l^d$ in the remainder of the paper for ease of analysis and presentation, but note that similar results are obtainable for other choices. 

Finally, for the SGCP model we assume a positive, continuous prior $p_{\lambda*}$ for $\lambda^*$ on $[0,\infty)$ satisfying \begin{equation}
\int_{\lambda'}^\infty p_{\lambda^*}(x)dx \leq C_0e^{-c_0(\lambda')^{\kappa}} \label{eq::lamstarassump}
\end{equation} for some constants $c_0,C_0,\kappa>0$ and all $\lambda'>0$. This condition is satisfied by, for instance, choosing a Gamma prior on $\lambda^*$.

\subsection{Additional Notation}

In the following section, we will derive results on the posterior distribution of $\lambda_0|{X}_{1:n}$ under the two models. We will denote the prior distributions as $\Pi(\cdot)$ and the posteriors as $\Pi(\cdot|X_{1:n})$. Certain results will be valid for the class of all continuous functions on $[0,1]^d$, which will be denoted $\mathcal{C}([0,1]^d)$, and others will hold for the class of all $\alpha$-H\"{o}lder continuous functions on $[0,1]^d$ denoted $\mathcal{C}^\alpha[0,1]^d$. 

Contraction results will inevitably depend on the particular filtering functions $\gamma_{1:n}$, therefore it is convenient to define versions of standard distances averaged with respect to $\gamma_{1:n}$. We have the averaged infity norm \begin{equation*}
    \Gamma_{n,\infty}(\lambda,\lambda')= \frac{1}{n}\sum_{i=1}^n ||\lambda\gamma_i - \lambda'\gamma_i||_\infty=\frac{1}{n}\sum_{i=1}^n \sup_{x \in [0,1]^d} |\lambda(x)\gamma_i(x)-\lambda'(x)\gamma_i(x)| ,
\end{equation*} averaged $L_2$ norm \begin{equation*}
    \Gamma_{n,2}(\lambda,\lambda') = \frac{1}{n}\sum_{i=1}^n ||\lambda\gamma_i - \lambda'\gamma_i||_2=\frac{1}{n}\sum_{i=1}^n \int_{[0,1]^d}(\lambda(x)\gamma_i(x)-\lambda'(x)\gamma_i(x))^2 dx,
\end{equation*} and square rooted averaged $L_2$ norm \begin{equation*}
    \Gamma_{n,2}^{1/2}(\lambda,\lambda') = \frac{1}{n}\sum_{i=1}^n ||\sqrt{\lambda\gamma_i} - \sqrt{\lambda'\gamma_i}||_2 =\frac{1}{n}\sum_{i=1}^n \int_{[0,1]^d}(\sqrt{\lambda(x)\gamma_i(x)}-\sqrt{\lambda'(x)\gamma_i(x)})^2 dx,
\end{equation*} for rate functions $\lambda, \lambda \in \mathcal{C}([0,1]^d)$. Using these definitions we can guarantee a rate of convergence appropriate to the level of filtering.

Finally let $N(\epsilon,\mathcal{S},l)$ denote the $\epsilon$-covering number of a set $\mathcal{S}$ with respect to distance $l$.

\section{Posterior Contraction Results} \label{sec::results}

In this section we state our results on the finite-time contraction of the posterior of the QGCP and SGCP models. Our results assert that given $n$ realisations of the NHPP, the expected posterior mass concentrated on functions outside a Hellinger-like ball of a given width will not exceed a transformation of the width of the ball. Theorem \ref{thm::inidQGCP} gives the result for the QGCP, and Theorem \ref{thm::inidSGCP} for the SGCP.

\begin{theorem} \label{thm::inidQGCP}
    	Suppose that $\lambda_0 \in \mathcal{C}^\alpha([0,1]^d)$ for some $\alpha>0$ and $\lambda_0:[0,1]^d \rightarrow [\lambda_{0,min},\infty)$. Suppose that the filtering functions $\gamma_{1:n}$ are known. Then for all sufficiently large $M,n>0$ the posterior under the QGCP satisfies \begin{align}
	&E_{\lambda_0}\bigg(\Pi\Big(\lambda: \frac{1}{n}\sum_{i=1}^n ||\sqrt{\lambda\gamma_i} - \sqrt{\lambda_0\gamma_i}||_2 \geq \sqrt{2}M\epsilon_n|{X}_{1:n}\Big)\bigg) = \tilde{o}\big(n^{\frac{-d}{4\alpha+d}}\big) \label{eq::inidQGCPconvergence}
	\end{align}
	for $\epsilon_n = 2\sqrt{||\lambda_0|_\infty} n^{-\alpha/(4\alpha+d)}(\log(n))^{\rho + d+1} + n^{-2\alpha/(4\alpha+d)}(\log(n))^{2\rho +2d+2} $ with  $\rho=\frac{1+d}{4+d/\alpha}$.
\end{theorem}

\begin{theorem} \label{thm::inidSGCP}
    Suppose that $\lambda_0 \in \mathcal{C}^\alpha([0,1]^d)$ for some $\alpha>0$ and $\lambda_0:[0,1]^d \rightarrow [\lambda_{0,min},\lambda_{0,max}]$. Suppose that the filtering functions $\gamma_{1:n}$ are known. Then for all sufficiently large $M,n>0$ the posterior under the SGCP satisfies \begin{align}
	&E_{\lambda_0}\bigg(\Pi\Big(\lambda: \frac{1}{n}\sum_{i=1}^n ||\sqrt{\lambda\gamma_i} - \sqrt{\lambda_0\gamma_i}||_2\geq \sqrt{2}M\epsilon_n|{X}_{1:n}\Big)\bigg) =\tilde{o}\big( n^{\frac{-d}{2\alpha+d}}\big) \label{eq::inidSGCPconvergence}
	\end{align}
	for $\epsilon_n = n^{-\alpha/(2\alpha+d)}(\log(n))^{\rho+d+1}$ with $\rho=\frac{1+d}{2+d/\alpha}$.
\end{theorem}

In each case analytical results free from ``little-$o$" notation and a specific value for the ``sufficiently large" conditions on $M$ and $n$ are given in the proofs in Sections \ref{sec::proofsketch} and \ref{sec::proofsketch2}.

The key difference between the two results is that for the QGCP we can only guarantee convergence on larger ball widths $\epsilon_n$ and at a slower rate $f_n$. Notice that under the QGCP the ball width is $\tilde{o}(n^{-\alpha/(4\alpha+d)})$ and the contraction rate is $\tilde{o}(n^{-d/(4\alpha+d)})$, whereas for the SGCP the ball width is $\tilde{o}(n^{-\alpha/(2\alpha+d)})$ and the contraction rate is $\tilde{o}(n^{-d/(2\alpha+d)})$. 

In the simplest setting where $\lambda_0 \in \mathcal{C}^1([0,1])$ - i.e. where we consider Lipschitz smooth functions on $d=1$ - this means we have a contraction rate of $\tilde{o}(n^{-1/5})$ on balls of width $\tilde{o}(n^{-1/5})$ for the QGCP and a contraction rate of $\tilde{o}(n^{-1/3})$ on balls of width $\tilde{o}(n^{-1/3})$ for the SGCP. The result on the SGCP is therefore tighter in two senses, we are able to say that the posterior mass shrinks quicker than for the QGCP and on the probability of being in a larger subspace (since the ball width $\epsilon_n$ is smaller, the area outside the ball is larger).

The different results arise as a consequence of the different transformation functions. For the posterior to contract at a given rate, we must demonstrate that the prior model satisfies certain properties related to this rate. Both models are built upon a GP $g$, and by considering the properties of $g$, we can verify that the SGCP and QGCP prior models meet the necessary conditions. 

The results of \citep{VanDerVaartVanZanten2009} demonstrate that for $g$ as described in Section \ref{sec::model}, relevant properties of $g$ can be shown, i.e. that the prior mass $g$ assigns to certain parts of the function space is bounded by sequences of a known form. It follows that appropriate transformations of these sequences can be used to show that the SGCP and QGCP priors also assign their prior mass across the function space in the required manner. The transformed sequences give rise to our ball widths $\epsilon_n$ which in turn influence the contraction rate. Since the SGCP and QGCP involve different transformations of $g$, we also require different transformations of the sequences for which desirable properties of $g$ hold, and therefore different results are obtained. 

More informally, the issue is that by applying a quadratic transformation to the GP over a logistic one, prior mass is dispersed more across the function space and the resulting posterior takes longer to contract around the true $\lambda_0$.

\subsection{Contraction of NHPP models under general priors}
Before we prove Theorems \ref{thm::inidQGCP} and \ref{thm::inidSGCP} we introduce a third result which gives a sufficient set of conditions on prior models to attain posterior contraction at a known rate under i.n.i.d. observations. Theorem \ref{thm::noniidgeneralPoisson} extends Theorem 1 of  \citep{GhosalVanDerVaart2007} to apply to for Poisson processes. The extension is in the same manner as the result of \citep{BelitserEtAl2013} extends Theorem 2 of \citep{GhosalVanDerVaart2001} for i.i.d. Poisson process realisations. In addition we retain the rate on the shrinkage of the posterior mass, as well as on the ball width, unlike these earlier papers.

\begin{theorem} \label{thm::noniidgeneralPoisson}
    Assume that $\lambda_0:[0,1]^d \rightarrow [\lambda_{0,min},\infty)$ and that filtering functions $\gamma_{1:n}$ are known. Suppose that for positive sequences $\delta,$  $\bar\delta_n \rightarrow 0$, such that $n\min(\delta_n,\bar\delta_n)^2 \rightarrow \infty$ as $n \rightarrow \infty $, it holds that there exist subsets $\Lambda_n \subset \mathcal{C}(S)$, some $n_0 \in \mathbb{N}$, and constants $c_1,c_2,c_3>0$, $c_4>1$, and $c_5>c_2+2$ such that \begin{align}
        \Pi_n \Big(\lambda: \Gamma_{n,\infty}(\lambda,\lambda_0) \leq \delta_n\Big) &\geq c_1 e^{-c_2n\delta_n^2} \label{eq::noniidPoissoncondition1} \\
        \sup_{\delta > \bar\delta_n} \log N \bigg(\frac{\delta}{36\sqrt{2}}, \sqrt{\Lambda_{n,\delta}},\Gamma_{n,2} \bigg) &\leq c_3n\bar\delta_n^2 \label{eq::noniidPoissoncondition2} \\
        \Pi_n(\Lambda \setminus \Lambda_n) &\leq c_4e^{-c_5n\delta_n^2} \label{eq::noniidPoissoncondition3}.
    \end{align}
    for all $n \geq n_0$ where $\Lambda_{n,\epsilon} = \Big\{\lambda \in \Lambda_n: h_n(p_\lambda,p_{\lambda_0})\leq \epsilon \Big\}$,   and $h_n(p_\lambda,p_{\lambda_0})$,  is given by \begin{displaymath}
h^2_n(p_\lambda,p_{\lambda'}) = \frac{1}{n}\sum_{i=1}^n 2\Bigg(1-E_{\lambda\gamma_i}\bigg( \sqrt{\frac{p(X^{(i)}|\lambda,\gamma_i)}{p(X^{(i)}|\lambda',\gamma_i)}}\bigg)\Bigg).
\end{displaymath} Then for $\epsilon_n = \max(\delta_n,\bar\delta_n)$ and any $C>0, J\geq 1, M\geq 2$, \begin{align}
        E_{\lambda_0}\bigg[\Pi_n\bigg(\lambda: \Gamma_{n,2}^{1/2}(\lambda,\lambda_0) \geq \sqrt{2}JM\epsilon_n|{X}_{1:n}\bigg)\bigg] &\leq  \frac{1}{C^2n\epsilon_n^2}+e^{-M^2n\epsilon_n^2/4} \nonumber\\
&\quad  + 2e^{-(M^2/2-c_3)n\epsilon_n^2} + \frac{2}{c_1}e^{-(c_2M^2J^2/4-C-1)n\epsilon_n^2} \label{eq::Belitserlikeconclusion}
    \end{align}
    for $n \geq \max(n_0,n_1,n_2,n_3)$ where $n_1 = \arg \min \{n: \epsilon_n \leq \lambda_{min} \}$, $n_2= \arg \min \{n: \epsilon_n \leq \frac{1}{\sqrt{2}M} \}$, and $n_3=\arg \min \{n: e^{-n\epsilon_n^2KM^2/4} \leq 1/2 \}$.
\end{theorem}

We prove this theorem in Section \ref{app::noniidPoissonproof}. This establishes that given the prior model satisfies certain conditions, the expected posterior mass assigned to rate functions outside an order $\epsilon_n$ width ball around $\lambda_0$ (measured with respect to an averaged $L_2$ distance) decreases at rate $o((n\epsilon_n^2)^{-1})$ for sufficiently large $n$. The conditions on the prior model are standard and are inherited from the conditions of Theorem 4 of \citep{GhosalVanDerVaart2007} required to show posterior contraction in a density estimation setting. Condition \eqref{eq::noniidPoissoncondition1}, the \emph{prior mass condition}, ensures that a sufficient proportion of the prior mass is assigned to functions close to $\lambda_0$. Condition \eqref{eq::noniidPoissoncondition2}, the \emph{entropy condition}, and condition \eqref{eq::noniidPoissoncondition3}, the \emph{remaining mass condition}, together prescribe that there exist subsets of the function space such that the entropy of these subsets is not too large, but the probability of lying outside these is also small.

Equipped with this general result, we are now in a position to prove Theorems \ref{thm::inidQGCP} and \ref{thm::inidSGCP} by demonstrating that the QGCP and SGCP models meet conditions \eqref{eq::noniidPoissoncondition1}, \eqref{eq::noniidPoissoncondition2}, and  \eqref{eq::noniidPoissoncondition3}.

\subsection{Proof of Theorem \ref{thm::inidQGCP}: Contraction of the QGCP model} \label{sec::proofsketch}

To prove Theorem \ref{thm::inidQGCP} we verify that the QGCP model  described in Section \ref{sec::model} meets the conditions of Theorem \ref{thm::noniidgeneralPoisson}. The following sections handles each condition in turn. Throughout we have \begin{align}
 \delta_n &=2\sqrt{||\lambda_0||_\infty} n^{-\alpha/(4\alpha+d)}\log^\rho(n) + n^{-2\alpha/(4\alpha+d)}\log^{2\rho}(n), \label{eq::deltaQGCP} \\
    \bar\delta_n &=2\sqrt{||\lambda_0||_\infty} n^{-\alpha/(4\alpha+d)}\log^{\rho+d+1}(n) + n^{-2\alpha/(4\alpha+d)}\log^{2\rho+2d+2}(n). \label{eq::deltabarQGCP}
\end{align}

\subsubsection{Prior Mass Condition}
The first condition, the so-called prior mass condition \eqref{eq::noniidPoissoncondition1} does not rely on the existence of particular subsets $\Lambda_n$, and can be verified by the following lemma, which we prove in Appendix  \ref{app::priormass}.

\begin{lemma} \label{lem::priormassbound}
    If $\lambda_0=g_0^2$ where $g_0 \in \mathcal{C}^\alpha([0,1]^d)$ for some $\alpha>0$ then under the QGCP model there exist constants $c_1,c_2>0$ for $\delta_n$ as defined in \eqref{eq::deltaQGCP} such that the prior satisfies \begin{equation*}
        \Pi(\lambda:||\lambda-\lambda_0||_\infty \leq \delta_n) \geq c_1e^{-c_2n\delta_n^2}
    \end{equation*}
     for all $n\geq 3$.
\end{lemma}

Then consider that since $\gamma_i \in [0,1]$ for all $i=1,...,n$, \begin{equation}
\Gamma_{n,\infty}(\lambda,\lambda_0)=\frac{1}{n}\sum_{i=1}^n ||\lambda\gamma_i - \lambda_0\gamma_i||_\infty \leq ||\lambda - \lambda_0||_\infty. \label{eq::inidtoiid}
\end{equation} Thus, by Lemma \ref{lem::priormassbound} we have that there exist constants $c_1,c_2>0$ such that \begin{equation*}
\Pi_n \Big(\lambda: \Gamma_{n,\infty}(\lambda,\lambda_0) \leq \delta_n\Big) \geq \Pi_n(\lambda: ||\lambda - \lambda_0||_\infty \leq \delta_n) \geq c_1 e^{-c_2n\delta_n^2},
\end{equation*} satisfying condition \eqref{eq::noniidPoissoncondition1}.

\subsubsection{Definition of Sieves}

We now define the subsets $\Lambda_n$ for which the QGCP satisfies the constraints of Theorem \ref{thm::noniidgeneralPoisson}. Let, \begin{equation}
\Lambda_n =(\mathcal{G}_n )^2 \label{eq::Lambdadef}
\end{equation} where \begin{equation}
    \mathcal{G}_n = \bigg[\beta_n \sqrt{\frac{\zeta_n}{\chi_n}}\mathbb{H}^{\zeta_n}_1 + \kappa_n \mathbb{B}_1\bigg] \cup  \bigg[ \bigcup_{a \leq \chi_n} (\beta_n\mathbb{H}_1^a)+\kappa_n\mathbb{B}_1 \bigg], \label{eq::Gdef}
\end{equation} $\mathbb{B}_1$ is the unit ball in $\mathcal{C}([0,1]^d)$ with respect to the uniform norm, and $\mathbb{H}_1^l$ is the unit ball of the RKHS $\mathbb{H}^l$ of the GP $g$ with covariance as given in \eqref{eq::covariance}. We define the sequences involved as follows, 
\begin{align*}
    \zeta_n &= L_2 n^{\frac{1}{d}\frac{2\alpha+d}{4\alpha+d}}(\log(n))^{2\rho/d} + L_3 n^{\frac{1}{d}\frac{\alpha+d}{4\alpha+d}}(\log(n))^{3\rho/d} + L_4 n^{\frac{1}{d}\frac{d}{4\alpha+d}}(\log(n))^{4\rho/d} \\
    \beta_n &= L_5 n^{\frac{1}{2}\frac{2\alpha+d}{4\alpha+d}}(\log(n))^{2\rho+\frac{d+1}{2}} + L_6 n^{\frac{1}{2}\frac{\alpha+d}{4\alpha+d}}(\log(n))^{3\rho+\frac{d+1}{2}} + L_7 n^{\frac{1}{2}\frac{d}{4\alpha+d}}(\log(n))^{4\rho+\frac{d+1}{2}}\\
    \kappa_n &= \frac{1}{3}\bar{\delta}_n, \quad \quad \quad
    \chi_n = \frac{\bar\delta_n}{6\tau \sqrt{d}\beta_n},
\end{align*} for constants \begin{align*}
L_2> (8c_5||\lambda_0||_\infty)/D_1, \quad \quad L_3>(8c_5\sqrt{||\lambda_0||_\infty})/D_1, \quad \quad L_4>2c_5/D_1
\end{align*}  such that $L_2+L_3+L_4>\max(A,e)$ and 
\begin{align*}
L_5 &\geq \max\bigg(\sqrt{\frac{16K_5L_2^d\mathcal{K}_1^{1+d}}{\log^{2\rho}(3)}},\sqrt{32||\lambda_0||_\infty c_5}, L_2^{1/3} \Big(\frac{8\max(1,\sqrt{||\lambda_0||_\infty})}{(3/36\sqrt{2})^{3/2}d^{1/4}\sqrt{2\tau}} \Big)^{2/3}\bigg) \\
L_6 &\geq \max\bigg(\sqrt{\frac{16K_5L_3^d\mathcal{K}_1^{1+d}}{\log^{3\rho}(3)}},\sqrt{32\sqrt{||\lambda_0||_\infty} c_5}\bigg), \enspace L_7 \geq \max\bigg(\sqrt{\frac{16K_5L_4^d\mathcal{K}_1^{1+d}}{\log^{4\rho}(3)}},\sqrt{8c_5}\bigg),
\end{align*} such that $L_5+L_6+L_7 > \frac{4L_1\max(1,\sqrt{||\lambda_0||_\infty})}{3\sqrt{||\mu||}}$, and $L_2L_5^3 > \bigg(\frac{8\max(1,||g_0||_\infty)}{(3/L_1)^{3/2}d^{1/4}\sqrt{2\tau}}\bigg)^2$ where $L_1=1/(36\sqrt{2})$ and $\mathcal{K}_1=\log\big(\frac{3(L_2+L_3+L_4)}{\min(1,2\sqrt{||\lambda_0|_\infty})} \big) + \frac{2\alpha d+ 2\alpha +d}{4\alpha d+d^2} + (4\rho-\rho/d-d-1)$.

The definition of $\mathcal{G}_n$ and these sequences is important as it allows general GP results of \citep{VanDerVaartVanZanten2009} to be applied. The extensive conditions on the constants are important to ensure that the results hold for finite values of $n$. 

\subsubsection{Entropy Condition}

The following lemma allows us to verify condition \eqref{eq::noniidPoissoncondition2} which stipulates that the log entropy of the subsets $\Lambda_n$ is not too large. The proof of this lemma is provided in Appendix \ref{app::entropy}. In particular it exploits an existing bound on the covering number of $\mathcal{G}_n$ with respect to the infinity norm from \citep{VanDerVaartVanZanten2009}. 

\begin{lemma}
    \label{lem::entropybound}
    For $\Lambda_n$ defined as in \eqref{eq::Lambdadef}, a constant $L_1>0$, and $\bar\delta_n$ as defined in \eqref{eq::deltabarQGCP}, there exists a constant $c_3>0$ such that \begin{equation*}
        \log N(L_1\bar\delta_n,\sqrt{\Lambda_n},||\cdot ||_2) \leq c_3 n \bar\delta_n^2,
    \end{equation*} for all $n$ such that \begin{align*}
      4||\lambda_0||_\infty \log^{2d+2-2\rho}(n) &\geq  \frac{m\sum_{i=2}^4 L_i^d}{2^{1+d}}\bigg(\log\Big(\frac{27\tau\sqrt{d}(\sum_{i=5}^7L_i)^3 \sum_{i=2}^4 L_i}{4||\lambda_0||_\infty^{3/2}}\Big)  +\Big(4+ \frac{12+d+d^2}{8\alpha d + 2d^2}\Big)\log(n)\bigg)^{1+d}
\end{align*} and\begin{equation*}
    2\log\bigg( \frac{6\sqrt{||\mu||}(L_5+L_6+L_7)}{2L_1\sqrt{||\lambda_0||_\infty}+L_1}\bigg) \leq 4||\lambda_0||_\infty n^{(4\alpha+2d)/(8\alpha+2d)}\log^{2\rho+2d+2}(n) - \log\bigg( n^{(6\alpha+d)/(4\alpha+d)}\log^{6\rho}(n) \bigg) .
\end{equation*}
\end{lemma}

To apply Lemma \ref{lem::entropybound}, notice that $\frac{1}{n}\sum_{i=1}^n ||\gamma_i \times \cdot ||_2 \leq  ||\cdot||_2$ since the functions $\gamma_i \in [0,1]$ for all $i=1,...,n$. It follows that \begin{align}
 N \bigg(\frac{\delta}{36\sqrt{2}}, \sqrt{\Lambda_{n,\delta}},\frac{1}{n}\sum_{i=1}^n ||\gamma_i \times \cdot||_2 \bigg) &\leq N \bigg( \frac{\delta}{36\sqrt{2}}, \sqrt{\Lambda_{n,\delta}}, ||\cdot||_2\bigg) \leq  N \bigg( \frac{\delta }{36\sqrt{2}}, \sqrt{\Lambda_n}, ||\cdot||_2 \bigg). \label{eq::entropyrelations}
\end{align} As any $\epsilon$-covering number is decreasing in $\epsilon$, it follows by Lemma \ref{lem::entropybound} that \begin{equation*}
     \sup_{\delta > \bar\delta_n} \log N \bigg(\frac{\delta}{36\sqrt{2}}, \sqrt{\Lambda_{n,\delta}},\frac{1}{n}\sum_{i=1}^n ||\gamma_i \times \cdot||_2 \bigg)  \leq c_3n\bar\delta_n^2.
     \end{equation*} Thus we have satisfied constraint \eqref{eq::noniidPoissoncondition2}.

\subsubsection{Remaining Mass Condition}

Finally, Lemma \ref{lem::remainingmassbound} below is sufficient to validate condition \eqref{eq::noniidPoissoncondition3} directly. Its proof is given in Appendix \ref{app::remainingmass}.

\begin{lemma}
    \label{lem::remainingmassbound}
    Under the QGCP model, with $\Lambda_n$ as defined in \eqref{eq::Lambdadef}, and $\delta_n$ as defined in \eqref{eq::deltaQGCP} there exist constants $c_4>0,c_5\geq c_2+4$ such that \begin{equation*}
        \Pi(\lambda: \lambda \notin \Lambda_n) \leq c_4e^{-c_5n\delta^2_n},
    \end{equation*} for all $n$ such that \begin{equation*}
        n^{(2\alpha+d)/(4\alpha +d)}\log^{2\rho}(n) \geq \frac{q_1}{4c_5||g_0||^2_\infty}\log\Big((L_2+L_3+L_4)n^{(2\alpha+d)/(4\alpha d + d^2)}\log^{4\rho/d}(n) \Big).
    \end{equation*}
\end{lemma}

\subsubsection{Conlcuding the Proof}
By Lemmas \ref{lem::priormassbound}, \ref{lem::entropybound}, and \ref{lem::remainingmassbound} and the definitions of $\delta_n$ and $\bar\delta_n$ therein we have that the conditions of Theorem \ref{thm::noniidgeneralPoisson} are satisfied. 
Thus, for the QGCP model \begin{align*}
E_{\lambda_0}\bigg[\Pi_n\bigg(\lambda: \Gamma_{n,2}^{1/2}(\lambda,\lambda_0) \geq \sqrt{2}JM\epsilon_n|\tilde{X}_{1:n}\bigg)\bigg] &\leq \frac{1}{C^2n\epsilon_n^2}+e^{-M^2n\epsilon_n^2/4} \nonumber\\
&\quad  + 2e^{-(M^2/2-c_3)n\epsilon_n^2} + \frac{2}{c_1}e^{-(c_2M^2J^2/4-C-1)n\epsilon_n^2}
\end{align*} holds with $\epsilon_n =\max(\delta_n,\bar\delta_n)=\bar\delta_n$ for any $C>0, J \geq 1, M\geq 2$. Specific values of the remaining constants can be extracted from Lemmas \ref{lem::priormassbound}, \ref{lem::entropybound}, \ref{lem::remainingmassbound}, and \ref{lem::modifiedvandervaart}. 

Then, so long as $M$ and $J$ are sufficiently large, the second, third and fourth terms on the RHS of equation \eqref{eq::Belitserlikeconclusion} decay much more quickly than the first and the bound is $\tilde{o}(n^{\frac{-d}{4\alpha+d}})$ as stated, for all $n$ such that the conditions of Theorem \ref{thm::noniidgeneralPoisson} and Lemmas \ref{lem::entropybound} and \ref{lem::remainingmassbound} are met. $\square$

\subsection{Proof of Theorem \ref{thm::inidSGCP}: Contraction of the SGCP model} \label{sec::proofsketch2}

Like the proof of Theorem \ref{thm::inidQGCP}, the proof of Theorem \ref{thm::inidSGCP} relies on demonstrating the the SGCP model described in Section \ref{sec::model} meets the conditions of Theorem \ref{thm::noniidgeneralPoisson}. In \citep{KirichenkoVanZanten2015} the conditions of Theorem 1 of \citep{BelitserEtAl2013} - the asymptotic and i.i.d. analogue of Theorem \ref{thm::noniidgeneralPoisson} - are verified for the SGCP model. However certain asymptotic arguments are used in said proof. In the following sections we handle each condition of Theorem \ref{thm::noniidgeneralPoisson} in turn under our setting. Throughout we have \begin{align}
\delta_n &= n^{-\alpha/(2\alpha+d)}(\log(n))^{(1+d)/(2+d/\alpha)} \label{eq::deltaSGCP} \\
\bar\delta_n &= n^{-\alpha/(2\alpha+d)}(\log(n))^{(1+d)/(2+d/\alpha)+d+1} \label{eq::deltabarSGCP}
\end{align}

\subsubsection{Prior Mass Condition}
For the SGCP model, the prior mass condition \eqref{eq::noniidPoissoncondition1} can be verified by the following lemma which we prove in Appendix \ref{app::priormassSGCP}. 

\begin{lemma} \label{lem::priormassSGCP}
If $\lambda_0=||\lambda_0||_\infty \sigma(g_0)$ where $g_0 \in \mathcal{C}^\alpha([0,1]^d)$ for some $\alpha>0$ then under the SGCP model there exist constants $c_1,c_2$ for $\delta_n$ as defined in \eqref{eq::deltaSGCP} such that the prior satisfies \begin{equation*}
        \Pi(\lambda:||\lambda-\lambda_0||_\infty \leq \delta_n) \geq c_1e^{-c_2n\delta_n^2}
    \end{equation*}
    for all $n\geq 3$.
\end{lemma}
\noindent
Then, by \eqref{eq::inidtoiid} and Lemma \ref{lem::priormassSGCP} we have that there exist constants $c_1,c_2>0$ such that \begin{equation*}
\Pi \bigg( \lambda: \Gamma_{n,\infty}(\lambda,\lambda_0) \leq \delta_n\bigg) \geq c_1e^{-c_2n\delta_n^2}
\end{equation*} and we have shown condition \eqref{eq::noniidPoissoncondition1} is satisfied under the SGCP model.

\subsubsection{Definition of Sieves}
We now define the sets $\Lambda_n$ such that the remaining coniditions hold. Consider, \begin{equation}
    \Lambda_n = \bigcup_{\lambda \leq \lambda_n} \lambda \sigma(\mathcal{G}_n) \label{eq::LambdaSGCP}
\end{equation} where \begin{equation}
    \mathcal{G}_n = \bigg[\beta_n \sqrt{\frac{\zeta_n}{\chi_n}}\mathbb{H}^{\zeta_n}_1 + \kappa_n \mathbb{B}_1\bigg] \cup  \bigg[ \bigcup_{a \leq \chi_n} (\beta_n\mathbb{H}_1^a)+\kappa_n\mathbb{B}_1 \bigg], 
\end{equation} $\mathbb{B}_1$ is the unit ball in $\mathcal{C}([0,1]^d)$ with respect to the uniform norm, and $\mathbb{H}_1^l$ is the unit ball of the RKHS $\mathbb{H}^l$ of the GP $g$ with covariance as given in \eqref{eq::covariance}. Though the structure of the sieves $\mathcal{G}_n$ is the same as in the proof of Theorem \ref{thm::inidQGCP}, the sequences are defined differently, as below
\begin{align*}
    \zeta_n &= L_8n^{\frac{1}{2\alpha+d}}(\log(n))^{2\rho/d}, \quad     \beta_n = L_9n^{\frac{d}{2(2\alpha+d)}}(\log(n))^{d+1+2\rho}, \\
    \lambda_n &= L_{10}n^{\frac{d}{\kappa(2\alpha+d)}}(\log(n))^{4\rho/\kappa}, \quad
    \kappa_n = \frac{1}{3}\bar\delta_n, \quad
    \chi_n = \frac{\kappa_n}{2\tau\sqrt{d}\beta_n},
\end{align*} for constants \begin{equation*}
L_8 > \max\Bigg( A, 1, \bigg(\frac{2c_5}{D_1} \bigg)^{1/d} \Bigg), \quad L_9 \geq \sqrt{8c_5}, \quad L_{10} > \bigg( \frac{c_5}{c_0}\bigg)^{1/\rho} 
\end{equation*} such that \begin{align*}
L_8L_9^3L_{10}^{3/2} > \frac{2}{(6cL_1)^{3/2}\tau\sqrt{d}}, \quad
L_9L_{10}^{1/2} > \frac{1}{6cL_1\sqrt{||\mu||}},
\end{align*} where $L_1 = 1/(36\sqrt{2})$, $c=2^{-5/2}$, and $\kappa$ is a positive constant.

\subsubsection{Entropy Condition}

The following lemma will allow us to verify condition \eqref{eq::noniidPoissoncondition2}. We prove it in Appendix \ref{app::entropySGCP}. \begin{lemma} \label{lem::entropySGCP}
For $\Lambda_n$ as defined in \eqref{eq::LambdaSGCP}, a constant $L_1>0$ and $\bar\delta_n$ as defined in \eqref{eq::deltabarSGCP}, there exists a constant $c_3>0$ such that \begin{equation*}
\log N(L_1\bar\delta_n, \sqrt{\Lambda_n},||\cdot||_2) \leq c_3n\bar\delta_n^2,
\end{equation*} for all $n$ such that \begin{align}
\log^{2d+2}(n) &>    K_1L_8^d \bigg(\log(\sqrt{2\tau L_8 L_9^3}L_{10}^{3/4}d^{1/4})+\frac{\kappa(6d+6\alpha+2)+3d}{4\kappa(2\alpha+d)}\log(n) \nonumber \\
&\quad \quad \quad \quad \quad \quad +\log\Big(\log^{3\rho/2 +3\rho/\kappa+\rho/d-d-1}(n)\Big)\bigg)^{1+d} \label{eq::SGCPentropy1}\\
n^{\frac{d}{2\alpha+d}} &> \max \Bigg( 2\log(12cL_1L_9L_{10}^{1/2}), 2\log(L_1L_{10}^{1/2}) \Bigg)+1. \label{eq::SGCPentropy2}
\end{align}
\end{lemma}

Then, as in the proof of Theorem \ref{thm::inidQGCP}, using \eqref{eq::entropyrelations} and Lemma \ref{lem::entropySGCP} we have \begin{equation*}
     \sup_{\delta > \bar\delta_n} \log N \bigg(\frac{\delta}{36\sqrt{2}}, \sqrt{\Lambda_{n,\delta}},\frac{1}{n}\sum_{i=1}^n ||\gamma_i \times \cdot||_2 \bigg)  \leq c_3n\bar\delta_n^2,
     \end{equation*} verifying condition \eqref{eq::noniidPoissoncondition2}.

\subsubsection{Remaining Mass Condition}

Finally, Lemma \ref{lem::remainingmassSGCP} below is sufficient to validate condition \eqref{eq::noniidPoissoncondition3} directly. Its proof is given in Appendix \ref{app:remainingmassSGCP}. 

\begin{lemma} \label{lem::remainingmassSGCP}
Under the SGCP model, with $\Lambda_n$ as defined in \eqref{eq::LambdaSGCP}, and $\delta_n$ as defined in \eqref{eq::deltabarSGCP} there exist constants $c_4>0$, $c_5 \geq c_2+4$ such that \begin{equation*}
\Pi(\lambda: \lambda \notin \Lambda_n) \leq c_4 e^{-c_5n\delta_n^2}
\end{equation*} for all $n$ such that  \begin{align}
\log^{2\rho}(n) &> \frac{16K_5D_1L_8^d}{L_9^2}\Bigg(\frac{\log(\sqrt{L_{10}}L_8)+\log\big(n^{\frac{2\alpha\kappa+ 2\kappa+d}{2\kappa(2\alpha+d)}}\log^{\rho(2/\kappa+2/d-1)-d-1}(n)\big) }{\log(n)} \Bigg)^{1+d}, \label{eq::SGCPncondition2} \\
n^{\frac{d}{2\alpha+d}}& >  \frac{1}{c_5}\big(\log(L_8^{q_1-d+1})+1 \big), \label{eq::SGCPncondition3}
\end{align}
\end{lemma}

\subsubsection{Concluding the Proof}
Thus, the three conditions \eqref{eq::noniidPoissoncondition1}, \eqref{eq::noniidPoissoncondition2}, and \eqref{eq::noniidPoissoncondition3} of Theorem \ref{thm::noniidgeneralPoisson} are satisfied by the SGCP model, 
and we have  \begin{align*}
E_{\lambda_0}\bigg[\Pi_n\bigg(\lambda: \Gamma_{n,2}^{1/2}(\lambda,\lambda_0) \geq \sqrt{2}JM\epsilon_n|\tilde{X}_{1:n}\bigg)\bigg] &\leq \frac{1}{C^k(n\epsilon_n^2)^{k/2}}+e^{-M^2n\epsilon_n^2/4} \nonumber\\
&\quad  + 2e^{-(M^2/2-c_3')n\epsilon_n^2} + \frac{2}{c_1'}e^{-(c_2'M^2J^2/4-C-1)n\epsilon_n^2}
\end{align*} with $\epsilon_n =\max(\delta_n,\bar\delta_n)=\bar\delta_n$.  Then, so long as $M$ and $J$ are sufficiently large, the second, third and fourth terms on the RHS of equation \eqref{eq::Belitserlikeconclusion} decay much more quickly than the first and the bound is $\tilde{o}(n^{\frac{-d}{2\alpha+d}})$ as stated, for all $n$ such that the conditions of Theorem \ref{thm::noniidgeneralPoisson} are met and that \eqref{eq::SGCPentropy1}, \eqref{eq::SGCPentropy2}, \eqref{eq::SGCPncondition2}, and \eqref{eq::SGCPncondition3} hold. $\square$ 

\subsection{Proof of Theorem \ref{thm::noniidgeneralPoisson}: Generic contraction in NHPPs} \label{app::noniidPoissonproof}
The proof of Theorem \ref{thm::noniidgeneralPoisson} depends on a general result for convergence of posterior parameter estimation given i.n.i.d. observations. Such a result is given in  \citep{GhosalVanDerVaart2007}, but without a finite time rate on the probability. We restate their result below as Theorem \ref{thm::generalinid} but with a rate included.

Consider as in \cite{GhosalVanDerVaart2007}, a model in which a parameter $\theta_0 \in \Theta$ gives rise to a i.n.i.d sequence of data. The data at time $i$ are drawn independently from the data at other times from a distribution $P_i^{\theta}$, which we assume admits a density $p_i^\theta$ with respect to a dominating measure.

We define the following subsets of the parameter space for $n\geq 1$ and $k>1$ \begin{equation*}
B_n(\theta_0,\epsilon;k)=\bigg\{\theta \in \Theta: \frac{1}{n}\sum_{i=1}^n K_i(\theta_0,\theta)\leq \epsilon^2, \enspace \frac{1}{n}\sum_{i=1}^n V_{k,0;i}(\theta_0,\theta) \leq C_k\epsilon^k \bigg\}
\end{equation*} where $K_i(\theta_0,\theta)=\int p_i^{\theta_0} \log(p_i^{\theta_0}/p_i^\theta)d\mu$ is the Kullback-Leibler divergence and $V_{k,0;i}(\theta_0,\theta)=\int p_i^{\theta_0}|\log(p_i^{\theta_0}/p_i^\theta)-K(\theta_0,\theta)|^k d\mu$ is a variance discrepancy measure. Furthermore, let $d_n$ be the averaged Hellinger distance, defined by \begin{displaymath}
d_n^2(\theta,\theta')=\frac{1}{n}\sum_{i=1}^n \int (\sqrt{p_{\theta,i}}-\sqrt{p_{\theta',i}})^2d\mu_i.
\end{displaymath} Our modified version of Theorem 4 of \citep{GhosalVanDerVaart2007} is as below.

\begin{theorem} \label{thm::generalinid}
    Suppose $Y_i \sim P_i^\theta$ independently for $i=1,...,n$ and let $d_n$ be defined as the average Hellinger distance. Further, suppose that for a sequence $\epsilon_n \rightarrow 0$ such that $n\epsilon_n^2$ is bounded away from 0, some $k >1$, all sufficiently large $j \in \mathbb{N}$, constants $c_1, c_2, c_3>0$, and sets $\Theta_n \subset \Theta$, the following conditions hold:

\begin{align}
\frac{\Pi_n(\theta \in \Theta_n: j\epsilon_n < d_n(\theta,\theta_0)\leq 2j\epsilon_n)}{\Pi_n(B^*_n(\theta_0,\epsilon_n;k))} &\leq c_1e^{\frac{c_2n\epsilon_n^2 j^2}{4}} \label{eq::ghosalinidpriormass}\\
\frac{\Pi_n(\Theta \setminus \Theta_n)}{\Pi_n(B_n^*(\theta_0,\epsilon_n;k))} &=o(e^{-2n\epsilon_n^2}) \label{eq::ghosalinidremainingmass} \\ 
\sup_{\epsilon>\epsilon_n} \log N \bigg(\frac{\epsilon}{36},\big\{\theta \in \Theta_n : d_n(\theta,\theta_0) < \epsilon \big\},d_n\bigg) &\leq  c_3n\epsilon_n^2.  \label{eq::ghosalinidentropy}
\end{align}
Then for any $C> 0, J\geq 1,$ and $M\geq 2$, \begin{align*}
\mathbb{E}_{\theta_0}\Pi_n(\theta: d_n(\theta,\theta_0) \geq JM\epsilon_n|Y^{(n)}) &\leq \frac{1}{C^k(n\epsilon_n^2)^{k/2}}+e^{-M^2n\epsilon_n^2/4} \\
&\quad \quad \quad + 2e^{-(M^2/2-c_3)n\epsilon_n^2} + \frac{2}{c_1}e^{-(c_2M^2J^2/4-C-1)n\epsilon_n^2}
\end{align*} for all $n$ such that $e^{-n\epsilon_n^2M^2/4} \leq 1/2$ .
\end{theorem}

The proof of Theorem \ref{thm::generalinid} is a modification of  proof of Theorem 4 of \citep{GhosalVanDerVaart2007}. We replace arguments that hold in the limit with finite-time versions and handle the introduction of the constants $c_1, c_2, c_3$, assumed to be 1 in \citep{GhosalVanDerVaart2007}. We can set the constant $K$ present in the original theorem to $1/2$ since we are dealing with the Hellinger distance.

\emph{Proof of Theorem \ref{thm::generalinid}:}
 By Lemmas 9 and 10 of \citep{GhosalVanDerVaart2007} and given conditions \eqref{eq::ghosalinidpriormass}, \eqref{eq::ghosalinidremainingmass}, and \eqref{eq::ghosalinidentropy}, we have for $n$ such that $e^{-n\epsilon_n^2M^2/4} \leq 1/2$ any $M \geq 2$, $J \geq 1$ and $C>0$, \begin{align*}
    &\mathbb{E}_{\theta_0}\Pi_n(\theta: d_n(\theta,\theta_0)\geq JM\epsilon_n|Y^{(n)}) \\
    &\leq \frac{1}{C^k(n\epsilon_n^2)^{k/2}}+e^{-M^2n\epsilon_n^2/4} + \frac{e^{-(M^2/2-c_3)n\epsilon_n^2}}{1-e^{-M^2n\epsilon_n^2/2}}  +\sum_{j \geq J} \frac{1}{c_1}e^{-n\epsilon_n^2(c_2M^2j^2/4 - C- 1)} \\
    &\leq \frac{1}{C^k(n\epsilon_n^2)^{k/2}}+e^{-M^2n\epsilon_n^2/4} + 2e^{-(M^2/2-c_3)n\epsilon_n^2}  +\frac{1}{c_1}e^{-n\epsilon_n^2(c_2M^2J^2/4 -C-1)}\sum_{j=0}^\infty {\big(e^{-n\epsilon_n^2(c_2M^2/4)}\big)}^{j^2} \\
    &\leq \frac{1}{C^k(n\epsilon_n^2)^{k/2}}+e^{-M^2n\epsilon_n^2/4} + 2e^{-(M^2/2-c_3)n\epsilon_n^2} +\frac{1}{c_1}e^{-n\epsilon_n^2(c_2M^2J^2/4 -C-1)}\sum_{j=0}^\infty {\big(e^{-n\epsilon_n^2(c_2M^2/4)}\big)}^{j^2} \\
    &\leq \frac{1}{C^k(n\epsilon_n^2)^{k/2}}+e^{-M^2n\epsilon_n^2/4}+ 2e^{-(M^2/2-c_3)n\epsilon_n^2} +\frac{e^{-n\epsilon_n^2(c_2M^2J^2/4 -C-1)}}{c_1(1-e^{-n\epsilon_n^2c_2M^2/4})} \\
    &\leq \frac{1}{C^k(n\epsilon_n^2)^{k/2}}+e^{-M^2n\epsilon_n^2/4}+ 2e^{-(M^2/2-c_3)n\epsilon_n^2} + \frac{2}{c_1}e^{-(c_2M^2J^2/4-C-1)n\epsilon_n^2}. \enspace \square
\end{align*} 

To apply Theorem \ref{thm::generalinid} we define averaged versions of the Hellinger distance, KL divergence and variance measure. Let $p_{\lambda\gamma_i}(N)=p(X^{(i)}|\lambda,\gamma_i)$ and we define the averaged Hellinger distance $h_n(p_\lambda,p_{\lambda'})$ by \begin{displaymath}
h^2_n(p_\lambda,p_{\lambda'}) = \frac{1}{n}\sum_{i=1}^n 2\Bigg(1-E_{\lambda\gamma_i}\bigg( \sqrt{\frac{p_{\lambda\gamma_i}(N)}{p_{\lambda'\gamma_i}(N)}}\bigg)\Bigg),
\end{displaymath}
the averaged KL-divergence as 
\begin{displaymath}
k_n(p_\lambda,p_{\lambda'}) = -\frac{1}{n}\sum_{i=1}^nE_{\lambda'\gamma_i}\bigg( \log\Big(\frac{p_{\lambda\gamma_i}(N)}{p_{\lambda'\gamma_i}(N)}\Big)\bigg),
\end{displaymath} and variance measure as \begin{displaymath}
v_n(p_\lambda,p_{\lambda'})=\frac{1}{n}\sum_{i=1}^n Var_{\lambda'\gamma_i}\bigg( \log\Big(\frac{p_{\lambda\gamma_i}(N)}{p_{\lambda'\gamma_i}(N)}\Big)\bigg).
\end{displaymath}

Through component-wise application of the relations in Section A.1 of \citep{BelitserEtAl2013} we have deterministic expressions for these quantities as \begin{align*}
h_n^2(p_{\lambda},p_{\lambda'}) &= \frac{1}{n}\sum_{i=1}^n 2\Bigg(1- \exp\bigg\{-\frac{1}{2}\int_{R_i}\Big(\sqrt{\lambda(t)\gamma_i(t)}-\sqrt{\lambda'(t)\gamma_i(t)}\Big)^2 dt \bigg\}\Bigg), \\
k_n(p_{\lambda},p_{\lambda'}) &= \frac{1}{n}\sum_{i=1}^n \Bigg(\int_{R_i}(\lambda(t)-\lambda'(t))\gamma_i(t) dt + \int_{R_i}\lambda'(t)\gamma_i(t)\log \Big(\frac{\lambda'(t)}{\lambda(t)}\Big)dt\Bigg), \\
v_n(p_\lambda,p_{\lambda'}) &= \frac{1}{n}\sum_{i=1}^n\int_{R_i}\lambda'(t)\gamma_i(t)\log^2\Big(\frac{\lambda'(t)}{\lambda(t)} \Big)dt.
\end{align*}
Lemma 1 of \citep{BelitserEtAl2013} gives bounds on the non-averaged versions of these quantities, but as the bounds will hold for each component of the average, we can trivially extend these results to give the following inequalities: \begin{align}
\frac{1}{\sqrt{2}n}\sum_{i=1}^n \big(||\sqrt{\lambda\gamma_i}-\sqrt{\lambda'\gamma_i}||_2 \wedge 1\big) \leq h_n(p_\lambda, p_{\lambda'}) &\leq \frac{\sqrt{2}}{n}\sum_{i=1}^n \big(||\sqrt{\lambda\gamma_i}-\sqrt{\lambda'\gamma_i}||_2 \wedge 1\big) \label{eq::lem1.1}\\
k_n(p_\lambda,p_{\lambda'}) &\leq \frac{3}{n}\sum_{i=1}^n  ||\sqrt{\lambda\gamma_i}-\sqrt{\lambda'\gamma_i}||^2_2 + v_n(p_\lambda,p_{\lambda'}) \label{eq::lem1.2}\\
\frac{1}{n}\sum_{i=1}^n ||\sqrt{\lambda\gamma_i}-\sqrt{\lambda'\gamma_i}||^2_2 &\leq \frac{1}{4n}\sum_{i=1}^n \int_{R_i} \gamma_i(s)(\lambda(s) \vee \lambda'(s))\log^2\Big(\frac{\lambda(s)}{\lambda'(s)} \Big) ds \label{eq::lem1.3}
\end{align}
where for numbers $x$ and $y$, the minimum is denoted $x \wedge y$ and the maximum is denoted $x \vee y$.

By assumption, $\lambda_0$ is bounded away from 0. It follows that any $\lambda \in \Lambda$ with $||\lambda_0-\lambda||_\infty \leq \lambda_{min}$ is also bounded away from 0, and that by the results \eqref{eq::lem1.2} and \eqref{eq::lem1.3} above $k_n(p_{\lambda_0}, p_{\lambda})$ and $v_n(p_{\lambda_0},p_\lambda)$ are both bounded by a constant times the averaged uniform norm $\frac{1}{n}\sum_{i=1}^n ||\lambda_0\gamma_i-\lambda\gamma_i||_\infty$. Therefore for $n \geq n_1$ the ball \begin{displaymath}
B_n^*(\epsilon_n) = \bigg\{\lambda \in \Lambda: k_n(p_{\lambda_0},p_{\lambda}) \leq \epsilon_n^2, v_n(p_{\lambda_0},p_{\lambda})\leq \epsilon_n^2 \bigg\}
\end{displaymath} is bounded by a multiple of the ball \begin{displaymath}
\bigg\{ \lambda \in \Lambda: \frac{1}{n}\sum_{i=1}^n ||\lambda_0 \gamma_i - \lambda \gamma_i||_\infty \leq \epsilon_n \bigg\}
\end{displaymath} for $\epsilon_n \leq \lambda_{min}$. It follows that for $n \geq n_1$, the condition \eqref{eq::noniidPoissoncondition1} implies \begin{equation}
    \Pi_n(B^*_n(\delta_n)) \geq c_1e^{-c_2n\delta_n^2}. \label{eq::inidKLballbound}
\end{equation}

By \eqref{eq::lem1.1} we have that \begin{displaymath}
N\bigg(\frac{\epsilon}{36},\Lambda_{n,\epsilon}, h_n \bigg) \leq N\bigg(\frac{\epsilon}{36\sqrt{2}},\sqrt{\Lambda_{n,\epsilon}},\frac{1}{n}\sum_{i=1}^n ||\cdot ||_2 \bigg)
\end{displaymath} where $\Lambda_{n,\epsilon} = \Big\{\lambda \in \Lambda_n: h_n(p_\lambda,p_{\lambda_0})\leq \epsilon \Big\}$. Thus the condition \eqref{eq::noniidPoissoncondition2} implies \eqref{eq::ghosalinidentropy}.

Combining these results we have \begin{displaymath}
\frac{\Pi_n(\lambda \in \Lambda_n: j \delta_n < h_n(p_\lambda,p_{\lambda_0}) \leq 2j\delta_n)}{\Pi_n(B^*_n(\delta_n))} \leq \frac{1}{\Pi_n(B^*_n(\delta_n))} \leq c_1e^{c_2n\delta_n^2}
\end{displaymath} by \eqref{eq::inidKLballbound} to satisfy condition \eqref{eq::ghosalinidpriormass}, and \begin{displaymath}
\frac{\Pi_n(\Lambda_n^c)}{\Pi_n(B^*_n(\delta_n))} \leq \frac{c_4e^{-c_5n\delta_n^2}}{c_1e^{-c_2n\delta_n^2}} = o(e^{-2n\delta_n^2})
\end{displaymath} by \eqref{eq::noniidPoissoncondition3} and \eqref{eq::inidKLballbound}
for $c_5-c_2 \geq 2$ to satisfy \eqref{eq::ghosalinidremainingmass}. Thus all the conditions of Theorem \ref{thm::generalinid} are satisfied by the assumptions of Theorem \ref{thm::noniidgeneralPoisson} and the conclusion of Theorem \ref{thm::generalinid} carries forward to Theorem \ref{thm::noniidgeneralPoisson} where we choose $k=2$. $\square$

\section{Conclusion}

We have derived finite time rates on the posterior contraction of the QGCP and SGCP models given i.n.i.d observations. This allows us to quantify the contraction of posterior estimates in the setting where events are not detected perfectly or the observation region is not sampled uniformly. As well as a new consistency result for the QGCP model, and the innovations of studying i.n.i.d data over i.i.d., the presentation of explicit rates on the posterior mass for the contraction of non-homogeneous Poisson process models is new. These results are of theoretical importance and practical interest in problems such as sequential decision making and experimental design.

We found that the SGCP model admitted a much tighter analysis than the QGCP model. For the simple setting of 1-dimensional 1-H\"{o}lder smooth rate functions, the SGCP model can be shown to have convergence of the near-optimal order $\tilde{o}(n^{-1/3}
)$. Our best result for the QGCP model only shows convergence of order $\tilde{o}(n^{-1/5})$. This discrepancy arises because of the different link functions used in the two models. In comparison to the bounded sigmoid function, the quadratic function induces a larger space of rate functions when the GP is transformed - meaning that wider sieves are required to give the desired results and the contraction guarantees are looser. Guarantees on the tightness of these bounds are currently unavailable, but this work provides some evidence to suggest that the SGCP model is superior to the QGCP in terms of rate of posterior contraction at least. This is an observation that would merit further empirical and analytical study in to the relationship between the models. We did not consider the LGCP model in this work as its exponential link function makes it very difficult to adapt the existing GP results of \citep{VanDerVaartVanZanten2009} into meaningful results in the NHPP posterior contraction setting. In particular, the high probability bound $\{||g-g_0||_\infty \leq \eta_{\beta_n}\}$ on the GP model, does not imply a useful bound on $||e^{g}-e^{g_0}||_\infty$ - the distance to be bounded in the prior mass condition for the LGCP - that gives useful contraction results.

We have focussed on particular choices of smoothness class, the link function used within the GCP construction and the width of the balls used in the contraction rate statements. There is of course potential to expand on these results by studying other choices. We believe however that the choices we have are consistent with the most common modelling choices in implementation of GCPs and useful for relating our results to the existing literature on posterior contraction of Bayesian nonparametric models.

\appendix

\section{Proof of Lemma \ref{lem::priormassbound}} \label{app::priormass}

Proving Lemma \ref{lem::priormassbound} relies on a bound on the uniform norm in the GP space. The following lemma gives a particular prior mass result which holds for all $n>0$ and uses a general term $\eta_{\beta,n}$ which fits with the analysis of both the SGCP and QGCP models. 
\begin{lemma}\label{lem::modifiedvandervaart}
	\label{lem::GPrate}
	If $g_0 \in \mathcal{C}^\alpha([0,1]^d)$ for some $\alpha>0$, then there exist constants $c_1, c_2>0$ such that  \begin{displaymath}
	\Pi\big(||g-g_0||_\infty \leq \eta_{\beta,n}\big) \geq c_1 e^{-c_2n\eta_{\beta,n}^\beta}
	\end{displaymath} for $\eta_{\beta,n}=n^{-\alpha/(\beta\alpha+d)}(\log(n))^{\rho_\beta}$ and $\rho_\beta=\frac{1+d}{\beta+d/\alpha}$ for all $n\geq 3$ where $\beta>1$.
\end{lemma}
\noindent
Furthermore, we rely on the following simple result which allows us to move between probabilistic bounds on the uniform norm of the GP and the squared GP.
\begin{lemma} \label{lem::simplenormresult}
    Let $w_1$ and $w_2$ be functions defined on $[0,1]^d$ such that $||w_2||_\infty$ is finite, and $c$ be a positive constant. Given the standard definition of the uniform norm, we have the following relation: \begin{equation*}
        \big\{||w_1-w_2||_\infty \leq c\big\} \Rightarrow \big\{||w_1^2 - w_2^2||_\infty \leq 2c||w_2||_\infty + c^2 \big\}.
    \end{equation*}
\end{lemma} 
\noindent We prove Lemma \ref{lem::modifiedvandervaart} in Appendix \ref{proof:modifiedvandervaart} and prove Lemma \ref{lem::simplenormresult} below.

\subsubsection*{Proof of Lemma \ref{lem::simplenormresult}:} We have: \begin{align*}
   \enspace ||w_1-w_2||_\infty &\leq c  && \\
   \Rightarrow w_1(x) &\leq w_2(x) + c \quad &&\forall \enspace x \in S \\
   \Rightarrow    w_1^2(x) &\leq w_2^2(x)+2cw_2(x) + c^2 \quad &&\forall \enspace x \in S \\
   \Rightarrow ||w_1^2 - w_2^2|| &\leq 2c||w_2||_\infty + c^2. && \square
\end{align*}

\subsubsection*{Proof of Lemma \ref{lem::priormassbound}:}

Recall the defintion $\delta_n = 2\eta_n ||g_0||_\infty + \eta_n^2$, with $\eta_n=\eta_{4,n}$. By definition we have: \begin{align*}
    &\quad \Pi\big(\lambda:||\lambda - \lambda_0||_{\infty} \leq \delta_n \big) \\
    &=  \Pi\big(g:||g^2-g_0^2||_\infty \leq 2\eta_n||g_0||_\infty + \eta_n^2 \big) \\
    &\geq \Pi\big( g: ||g-g_0||_\infty \leq \eta_n \big) \\
    &\geq c_1e^{-c_2n\eta_n^4} \geq c_1e^{-c_2n\delta_n^2},
\end{align*}  Here, the first inequality is due to Lemma \ref{lem::simplenormresult}. The second is by application of Lemma \ref{lem::modifiedvandervaart} and the third is by definition of $\delta_n$. $\square$

\section{Proof of Lemma \ref{lem::modifiedvandervaart}} \label{proof:modifiedvandervaart}

We will utilise the following result from Section 5.1 of \citep{VanDerVaartVanZanten2009}, which holds for a constant $H$ depending only on $g_0$ and $\mu$, a constant $K_2$ depending only on $g_0, \mu, \alpha, d$ and $D_1$ and any $\epsilon>0$ \begin{displaymath}
\Pi\big(||g- g_0||_\infty \leq 2\epsilon \big) \geq C_1 \exp\bigg\{-K_2 \bigg(\frac{1}{\epsilon}\bigg)^{d/\alpha}\Bigg(\log\bigg(\frac{1}{\epsilon}\bigg) \Bigg)^{1+d} \bigg\}\bigg(\frac{H}{\epsilon}\bigg)^{\frac{q_1+1}{\alpha}}.
\end{displaymath} Recall that $C_1$, $D_1$, and $q_1$ are constants from the assumption \eqref{eq::lengthscalecondition} on the length scale of the GP prior. 

Substituting the particular form of $\epsilon=\epsilon_n=n^{-\alpha/(\beta\alpha+d)}(\log(n))^{\rho}$ and $\rho=\frac{1+d}{\beta+d/\alpha}$ from Lemma \ref{lem::modifiedvandervaart} into the above we have:
\begin{align*}
\Pi\big(||g-g_0||_\infty \leq \epsilon_n\big) &\geq C_1\exp\Bigg\{ -2^{\frac{d}{\alpha}} K_2 n^{\frac{d}{\beta\alpha+d}}(\log(n))^{-\frac{d}{\alpha}\frac{1+d}{\beta+d/\alpha}} \bigg(\log\bigg(2n^{\frac{\alpha}{\beta\alpha +d}}(\log(n))^{-\frac{1+d}{\beta+d/\alpha}}\bigg)\bigg)^{1+d}\Bigg\} \\
&\quad \times \Bigg(2Hn^{\frac{\alpha}{\beta\alpha +d}}(\log(n))^{-\frac{1+d}{\beta+d/\alpha}}\Bigg)^{\frac{q_1+1}{\alpha}}, \\
\intertext{defining $Z(n)= \Big(Hn^{\frac{\alpha}{\beta\alpha +d}}(\log(n))^{-\frac{1+d}{\beta+d/\alpha}}\Big)^{\frac{q_1+1}{\alpha}}$ and expanding the logarithm,}
&= C_1 Z(n) \exp \Bigg\{-2^{\frac{d}{\alpha}}K_2 n^{\frac{d}{\beta\alpha+d}}(\log(n))^{-\frac{d}{\alpha}\frac{1+d}{\beta+d/\alpha}}  \bigg(\frac{\alpha}{\beta\alpha+d}\log(2n)- (\log(n))^{\frac{1+d}{\beta+d/\alpha}} \bigg)^{1+d} \Bigg\} \\
&\geq C_1 Z(n) \exp \Bigg\{-2^{\frac{d}{\alpha}}K_2 n^{\frac{d}{\beta\alpha+d}}(\log(n))^{-\frac{d}{\alpha}\frac{1+d}{\beta+d/\alpha}} \bigg(\frac{\alpha}{\beta\alpha+d}\log(2n)\bigg)^{1+d} \Bigg\} \\
\intertext{using $2\log(n) \geq \log(2n)$ for $n \geq 2$} &\geq C_1 Z(n) \exp \Bigg\{-2^{1+d/\alpha}K_2 n\cdot n^{-\frac{\beta\alpha}{\beta\alpha+d}}(\log(n))^{(1+d)\frac{\beta+d/\alpha -d/\alpha}{\beta+d/\alpha}}  \Bigg\}  \\
&=C_1Z(n)\exp\Bigg\{-2^{1+d/\alpha}K_2n\epsilon_n^\beta\Bigg\} \\
\intertext{letting $K_3 = \min_{n \geq 3}(Z(n))$}
&\geq C_1K_3\exp\Bigg\{-2^{1+d/\alpha} K_2 n\epsilon_n^\beta\Bigg\} = c_1e^{-c_2n\epsilon_n^\beta}, 
\end{align*} where $c_1=C_1K_3$, and $c_2=2^{1+d/\alpha}K_2$. $\square$

\section{Proof of Lemma \ref{lem::entropybound}} \label{app::entropy}

As $\sqrt{\Lambda_n}=\mathcal{G}_n$, the covering numbers $N(L_1\bar\delta_n, \sqrt{\Lambda_n},||\cdot||_2)$ and $N(L_1\bar\delta_n, \mathcal{G}_n, ||\cdot ||_2)$ are equivalent. It follows that \begin{equation*}
N(L_1\bar\delta_n, \sqrt{\Lambda_n},||\cdot||_2) \leq N(L_1\bar\delta_n, \mathcal{G}_n, ||\cdot ||_\infty).
\end{equation*}  Defining $\mathcal{G}_n$ as in \eqref{eq::Gdef} allows us to use the following result, (5.4) of \citep{VanDerVaartVanZanten2009}: \begin{equation*}
    \log N (L_1\bar\delta_n, \mathcal{G}_n,||\cdot||_\infty) \leq m\zeta_n^d\bigg(\log \frac{3^{3/2}d^{1/4}\beta_n^{3/2}\sqrt{2\tau \zeta_n}}{(L_1\bar\delta_n)^{3/2}} \bigg)^{1+d}+ 2\log \frac{6\beta_n\sqrt{||\mu||}}{L_1\bar\delta_n}
\end{equation*} for $||\mu||$ the total mass of the spectral measure $\mu$, $\tau^2$ as the second moment of $\mu$, positive constant $m$ depending only on $\mu$ and $d$, and given \begin{equation*}
    (3/L_1)^{3/2}d^{1/4}\beta_n^{3/2}\sqrt{2\tau \zeta_n} > 2\bar\delta_n^{3/2}, \quad \quad (3/L_1)\beta_n\sqrt{||\mu||} > \bar\delta_n.
\end{equation*}
By the definitions of $\beta_n$, and $\zeta_n$ we have that \begin{align*}
    m\zeta_n^d\bigg(\log \frac{3^{3/2}d^{1/4}\beta_n^{3/2}\sqrt{2\tau \zeta_n}}{(L_1\bar\delta_n)^{3/2}} \bigg)^{1+d} &\leq n\bar\delta_n^2 \\
    2\log \frac{6\beta_n\sqrt{||\mu||}}{L_1\bar\delta_n} &\leq n\bar\delta_n^2,
\end{align*} for the values of $n$ specified in the statement of  Lemma \ref{lem::entropybound}. It follows that the lemma is satisfied with $c_3=2$. $\square$

\section{Proof of Lemma \ref{lem::remainingmassbound}}  \label{app::remainingmass}

Firstly note that $\Pi(\lambda \notin \Lambda_n)=\Pi(g \notin \mathcal{G}_n)$. By a simplification of (5.3) of \citep{VanDerVaartVanZanten2009} to account for our assumption that $q_2=0$, we have \begin{equation*}
    \Pi(g \notin \mathcal{G}_n) \leq C_1\zeta_n^{q_1-d+1}e^{-D_1\zeta_n^d }+e^{-\beta_n^2/8}
\end{equation*} $\bar\delta_n <\delta_0$ for small $\delta_0>0$, and $\beta_n$, $\zeta_n$, and $\bar\delta_n$ satisfying \begin{equation*}
    \beta_n^2 > 16K_5\zeta_n^d\bigg(\log\Big(\frac{3\zeta_n}{\bar\delta_n}\Big)\bigg)^{1+d}, \quad \quad \zeta_n > 1, 
\end{equation*} for a constant $K_5$ depending only on $\mu$ and $g$. 
The definitions of $\beta_n$, $\delta_n$ and $\zeta_n$ give us the following relations, for a constant $c_5 = c_2+4$ \begin{align*}
    D_1 \zeta_n^d &\geq 2c_5n\delta_n^2, \enspace   \beta_n^2                 \geq 8c_5n\delta_n^2, \enspace
    \zeta_n^{q_1-d+1} \leq e^{c_5n\delta_n^2}, 
\end{align*} with the final of these holding for values of $n$ as specified in the statement of Lemma \ref{lem::remainingmassbound}. Using these we can obtain the necessary result as follows: \begin{align*}
    \Pi(\lambda:\lambda \notin \Lambda_n) =\Pi(g : g \notin \mathcal{G}_n) &\leq C_1\zeta_n^{q_1-d+1}e^{-D_1\zeta_n^d }+e^{-\beta_n^2/8} \\
                                        &\leq C_1e^{c_5n\delta_n^2}e^{-2c_5n\delta_n^2} + e^{-c_5n\delta_n^2} \\
                                        &= \big(C_1 + 1 \big)e^{-2c_5n\delta_n^2} \\
                                        &\leq c_4e^{-(c_2+4)n\delta_n^2}
\end{align*} for $c_4 = C_1+1$. $\square$

\section{Proof of Lemma \ref{lem::priormassSGCP}} \label{app::priormassSGCP}
Under the SGCP model we have \begin{align*}
     \Pi\bigg((\lambda^*,g):||\lambda^* \sigma(g) - \lambda_0 ||_\infty \leq \delta_n\bigg) \geq \Pi\bigg(\lambda^*:|\lambda^* - 2||\lambda_0||_\infty | \leq \frac{\delta_n}{2}\bigg)\Pi\bigg(g:||\sigma(g)-\sigma(g_0)||_\infty \leq \frac{\delta_n}{4||\lambda_0||_\infty}\bigg).
\end{align*} By the assumption that $\lambda^*$ has a positive continuous density, the first term on the RHS of the inequality can be bounded below by a constant times $\delta_n$, which can itself be lower bounded by a constant for finite $n$. The second term can be bounded below by $\Pi_n(||g-g_0||_\infty \leq \delta_n/(16||\lambda_0|_\infty))$ since $1/4$ is the Lipschitz constant of the sigmoid transformation.  Thus, by Lemma \ref{lem::modifiedvandervaart} (given in Appendix \ref{app::priormass}) we have: \begin{displaymath}
\Pi \Big(\lambda: \Gamma_{n,\infty}(\lambda,\lambda_0) \leq \delta_n\Big) \geq c_0'\Pi\bigg(g: ||g-g_0||_\infty \leq \frac{\delta_n}{16||\lambda_0||_\infty}\bigg) \geq c_1' e^{-nc_2'\delta_n^2}
\end{displaymath} for positive constants $c_1',c_2'$, showing condition \eqref{eq::noniidPoissoncondition1} is satisfied under the SGCP model.

\section{Proof of Lemma \ref{lem::entropySGCP}} \label{app::entropySGCP}

Define $\psi_n = \bar\delta_n/(2L_1\sqrt{\lambda_n})$. We have \begin{align}
\log N(L_1\bar\delta_n,\sqrt{\Lambda_n},||\cdot||_2) &= \log N(2\psi_n \sqrt{\lambda_n}, \sqrt{\Lambda_n},||\cdot||_2) \nonumber \\
&\leq \log N(\psi_n\sqrt{\lambda_n},[0,\lambda_n],\sqrt{|\cdot |}) + \log N(\psi_n/c,\mathcal{G}_n,||\cdot ||_\infty) \nonumber \\
&\leq \log\frac{1}{\psi_n} + \log N(\psi_n/c,\mathcal{G}_n,||\cdot ||_\infty) \label{eq::SGCPentropytool1}
\end{align} for $c=2^{-5/2}$, the Lipschitz constant of $\sqrt{\sigma}$. 

Then, as in the proof of Lemma \ref{lem::entropybound}, by equation (5.4) of \citep{VanDerVaartVanZanten2009}, we have for $B_n>0$, \begin{equation}
    \log N (B_n\bar\delta_n, \mathcal{G}_n,||\cdot||_\infty) \leq m\zeta_n^d\bigg(\log \frac{3^{3/2}d^{1/4}\beta_n^{3/2}\sqrt{2\tau \zeta_n}}{(B_n\bar\delta_n)^{3/2}} \bigg)^{1+d}+ 2\log \frac{6\beta_n\sqrt{||\mu||}}{B_n\bar\delta_n} \label{eq::SGCPentropytool2}
\end{equation}
 subject to the conditions \begin{align*}
    (3/B_n)^{3/2}d^{1/4}\beta_n^{3/2}\sqrt{2\tau \zeta_n} > 2\bar\delta_n^{3/2}, \quad (3/B_n)\beta_n\sqrt{||\mu||} > \bar\delta_n, \quad \zeta_n >A
\end{align*} for a constant $A>0$. These conditions hold by defintion for $n$ as specified by \eqref{eq::SGCPentropy1}, with $B_n=(2cL_1\sqrt{\lambda_n})^{-1}$. Then, combining \eqref{eq::SGCPentropytool1} and \eqref{eq::SGCPentropytool2} we have  \begin{align*}
\log N \big(L_1\bar\delta_n, \sqrt{\Lambda_{n}},||\cdot||_2 \big) &\leq \log  \frac{2L_1\sqrt{\lambda_n}}{\bar\delta_n} + m\zeta_n^d\bigg(\log \frac{(6cL_1)^{3/2}d^{1/4}\beta_n^{3/2}\sqrt{2\tau \lambda_n \zeta_n}}{\bar\delta_n^{3/2}} \bigg)^{1+d}\\
&\quad \quad + 2\log \frac{12cL_1\beta_n\sqrt{\lambda_n||\mu||}}{\bar\delta_n}.
\end{align*} For $n$ as specified by \eqref{eq::SGCPentropy2}, we have \begin{align*}
	\log  \frac{2L_1\sqrt{\lambda_n}}{\bar\delta_n} &< n \bar\delta_n^2,    \\
    m\zeta_n^d\bigg(\log \frac{(6cL_1)^{3/2}d^{1/4}\beta_n^{3/2}\sqrt{2\tau \lambda_n \zeta_n}}{\bar\delta_n^{3/2}} \bigg)^{1+d} &< n \bar\delta_n^2,\\
         2\log \frac{12cL_1\beta_n\sqrt{\lambda_n||\mu||}}{\bar\delta_n} &< n \bar\delta_n^2. 
\end{align*} Thus for $n$ satisfying \eqref{eq::SGCPentropy1} and \eqref{eq::SGCPentropy2} we have \begin{equation*}
\log N \bigg(L_1\bar\delta_n, \sqrt{\Lambda_{n}},||\cdot||_2 \bigg)\leq 3n\bar\delta_n^2
\end{equation*} proving Lemma \ref{lem::entropySGCP} with $c_3=3$. $\square$

\section{Proof of Lemma \ref{lem::remainingmassSGCP}} \label{app:remainingmassSGCP}

As in \citep{KirichenkoVanZanten2015}, we may decompose the probability of interest \begin{align*}
\Pi\Big(\lambda: \lambda \notin \Lambda_n\Big) &= \Pi\Big((\lambda^*,g):\lambda^* \sigma(g) \notin \Lambda_n\Big) \\
&\leq \int_0^{\lambda_n} \Pi\Big((\lambda^*,g):\lambda^* \sigma(g) \notin \Lambda_n\Big)p_{\lambda^*}(\lambda)d\lambda + \int_{\lambda_n}^{\infty} p_{\lambda^*}(\lambda) d\lambda \\
&\leq \Pi\Big(g:g \notin \mathcal{G}_n\Big) + C_0 e^{-c_0 \lambda_n^\rho},
\end{align*} by the assumption \eqref{eq::lamstarassump}.  As utilised in the proof of Lemma \ref{lem::remainingmassbound}, equation (5.3) of \citep{VanDerVaartVanZanten2009} states that \begin{equation*}\Pi(g \notin \mathcal{G}_n) \leq C_1\zeta_n^{q_1-d+1}e^{-D_1\zeta_n^d }+e^{-\beta_n^2/8}
\end{equation*} given conditions \begin{equation*}
    \beta_n^2 > 16 K_5 \zeta_n^d \bigg(\log \Big(\frac{\lambda_n^{1/2} \zeta_n}{\bar\delta_n} \Big) \bigg)^{1+d}, \quad \zeta_n > 1
\end{equation*}  which are satisfied by our earlier definitions, for a constant $K_5$ depending only on $\mu$ and $g$. Then for $n$ as specified by equations \eqref{eq::SGCPncondition2} and \eqref{eq::SGCPncondition3}, we have the following results \begin{align*}
    c_0\lambda_n^\rho > c_5n\delta_n^2, \quad  D_1\zeta_n^d \geq 2c_5n\delta_n^2, \quad   \zeta_n^{q_1-d+1} \leq e^{c_5n\delta_n^2}, \quad   \beta_n^2 \geq 8c_5n\delta_n^2.
\end{align*} The required result then follows. $\square$

\section{Verifying conditions on sieves}
Throughout the analysis the sequences used in defining the sieves are subject to numerous conditions and assumptions, in order that we may demonstrate the conditions of Theorem \ref{thm::noniidgeneralPoisson} are met for the GCP models. By choosing $L_{2:10}$ as specified in the main body, these conditions are met by definition for values of $n$ as specified. There are numerous such conditions to verify, and doing so can be non-trivial. In this section we show the link between the conditions and constraints on $L_{2:10},n$ and demonstrate fully that the necessary results hold.

\subsection{QGCP model}
Recall, the definitions of the following sequences: \begin{align*}
    \delta_n &=2||g_0||_\infty n^{-\alpha/(4\alpha+d)}\log^\rho(n) + n^{-2\alpha/(4\alpha+d)}\log^{2\rho}(n), \\
    \bar\delta_n &=2||g_0||_\infty n^{-\alpha/(4\alpha+d)}\log^{\rho+d+1}(n) + n^{-2\alpha/(4\alpha+d)}\log^{2\rho+2d+2}(n) \\
    \zeta_n &= L_2n^{(2\alpha+d)/(4\alpha d+d^2)}\log^{2\rho/d}(n) + L_3n^{(\alpha d+d^2)/(4\alpha d+d^2)}\log^{3\rho/d}(n) +L_4n^{d/(4\alpha+d)}\log^{4\rho/d}(n) \\
    \beta_n &= L_5n^{(2\alpha+d)/(8\alpha+2d)}\log^{2\rho+(d+1)/2}(n) + L_6n^{(\alpha+d)/(8\alpha+2d)}\log^{3\rho+(d+1)/2}(n) \\
    &\quad \quad \quad +L_7n^{d/(8\alpha+2d)}\log^{4\rho+(d+1)/2}(n)
\end{align*} with $L_2,...,L_7$ satisfying \begin{align*}
    L_2+L_3+L_4 &> \max(A,e) \\
    L_2L_5^3 &> \bigg(\frac{8\max(1,||g_0||_\infty)}{(3/L_1)^{3/2}d^{1/4}\sqrt{2\tau}}\bigg)^2 \\
    L_5+L_6+L_7& > \frac{4L_1\max(1,||g_0||_\infty)}{3\sqrt{||\mu||}} \\
    L_2 &\geq (8c_5||g_0||_\infty^2)/D_1 \\
    L_3 &\geq (8c_5||g_0||_\infty)/D_1 \\
    L_4 &\geq 2c_5/D_1 \\
    L_5 &\geq \max\bigg(\sqrt{\frac{16K_5L_2^d\mathcal{K}_1^{1+d}}{\log^{2\rho}(3)}},\sqrt{32||g_0||_\infty^2 c_5}\bigg) \\
    L_6 &\geq \max\bigg(\sqrt{\frac{16K_5L_3^d\mathcal{K}_1^{1+d}}{\log^{3\rho}(3)}},\sqrt{32||g_0||_\infty c_5}\bigg) \\
    L_7 &\geq \max\bigg(\sqrt{\frac{16K_5L_4^d\mathcal{K}_1^{1+d}}{\log^{4\rho}(3)}},\sqrt{8c_5}\bigg)
\end{align*} for $n \geq \max(3,n_3,n_4,n_5)$. Here $n_3$ is the smallest integer $n$ such that\begin{align*}
      4||g_0||_\infty^2 \log^{2d+2-2\rho}(n) &\geq  \frac{m\sum_{i=2}^4 L_i^d}{2^{1+d}}\bigg(\log\Big(\frac{27\tau\sqrt{d}(\sum_{i=5}^7L_i)^3 \sum_{i=2}^4 L_i}{4||g_0||_\infty^3}\Big) \\
    &\quad \quad \quad \quad +\Big(4+ \frac{12+d+d^2}{8\alpha d + 2d^2}\Big)\log(n)\bigg)^{1+d}
\end{align*} $n_4$ is the smallest integer $n$ such that \begin{equation*}
     2\log\bigg( \frac{6\sqrt{||\mu||}(L_5+L_6+L_7)}{2L_1||g_0||_\infty+L_1}\bigg) \leq 4||g_0||_\infty^2 n^{(4\alpha+2d)/(8\alpha+2d)}\log^{2\rho+2d+2}(n) - \log\bigg( n^{(6\alpha+d)/(4\alpha+d)}\log^{6\rho}(n) \bigg) 
\end{equation*} and $n_5$ is the smallest integer $n$ such that \begin{equation*}
     n^{(2\alpha+d)/(4\alpha +d)}\log^{2\rho}(n) \geq \frac{q_1}{4c_5||g_0||^2_\infty}\log\Big((L_2+L_3+L_4)n^{(2\alpha+d)/(4\alpha d + d^2)}\log^{4\rho/d}(n) \Big).
\end{equation*}

In the remainder of this subsection, we show that the following conditions, which are all restatements of required results in our main analysis, hold for the sequences described above. \begin{align}
    \zeta_n &> \max(A,1) \label{QGCP1}\\
    (3/L_1)^{3/2}d^{1/4}\beta_n^{3/2}\sqrt{2\tau\zeta_n} &> 2\bar\delta_n^{3/2} \label{QGCP2}\\
    (3/L_1)\beta_n\sqrt{||\mu||} &> \bar\delta_n \label{QGCP3}\\
    m\zeta_n^d \Bigg(\log\bigg(\frac{(3/L_1)^{3/2}d^{1/4}\beta_n^{3/2}\sqrt{2\tau\zeta_n}}{\bar\delta_n^{3/2}} \bigg) \Bigg)^{1+d} &\leq n\bar\delta_n^2 \label{QGCP4}\\
    2\log\bigg( \frac{6\beta_n\sqrt{||\mu||}}{L_1\bar\delta_n}\bigg) &\leq n\bar\delta_n^2 \label{QGCP5}\\
    \beta_n^2 &> 16K_5\zeta_n^d\Bigg(\log \bigg(\frac{3\zeta_n}{\bar\delta_n} \bigg) \Bigg)^{1+d} \label{QGCP6}\\
    D_1\zeta_n^d\bigg(\log^{q_2}(\zeta_n)\bigg) &\geq 2c_5n\delta_n^2 \label{QGCP7}\\
    \beta_n^2 &\geq 8c_5n\delta_n^2 \label{QGCP8}\\
    \zeta_n^{q_1-d+1} &\leq \exp(c_5n\delta_n^2), \label{QGCP9}
\end{align} 

\subsubsection{Verifying \eqref{QGCP1}}
For $n=3$, $\log(n)>1$ thus $\zeta_n > L_2 +L_3 + L_4$ for all $\alpha \in [0,1]$, and $d\geq 1$. It follows that \eqref{QGCP1} is satisfied for $n=3$ given $L_2+L_3+L_4 > \max(A,e)$. To show it holds for all $n>3$ we simply note that $\zeta_n$ is an increasing function.
\subsubsection{Verifying \eqref{QGCP2}}
First consider, \begin{align*}
    \bar\delta_n &= 2||g_0||_\infty n^{-\alpha/(4\alpha+d)}\log^{\rho+d+1}(n) + n^{-2\alpha/(4\alpha+d)}\log^{2\rho+2d+2}(n) \\
    &\leq 4\max(1,||g_0||_\infty)n^{-\alpha/(4\alpha+d)}\log^{2\rho+2d+2}(n) \\
    \Rightarrow \quad 2\bar\delta_n^{3/2} &\leq (4\max(1,||g_0||_\infty))^{3/2} n ^{-3\alpha/(8\alpha+2d)} \log^{3\rho+3d+3}(n))^{3/2} \\
       &= (4\max(1,||g_0||_\infty))^{3/2}n^{-3\alpha/(8\alpha +2d)}\log^{3\rho+3(d+1)/4}(n)\log^{9(d+1)/4}(n) \\
\end{align*}
Let $z_1=(3/L_1)^{3/2}d^{1/4}\sqrt{2\tau}$, \begin{align*}
   z_1 \sqrt{\beta_n^3\zeta_n} &\geq z_1\log^{4\rho+3(d+1)/4}(n)\sqrt{\big(L_5n^{\frac{2\alpha+d}{8\alpha+2d}} + L_6n^{\frac{\alpha+d}{8\alpha+2d}}+L_7n^{\frac{d}{8\alpha+2d}}\big)^3\big(L_2n^{\frac{2\alpha+d}{4\alpha d+d^2}} + L_3n^{\frac{\alpha+d}{4\alpha d+d^2}} +L_4n^{\frac{d}{4\alpha d+d^2}}\big)} \\
   &\geq z_1\log^{3\rho+3(d+1)/4}(n)\sqrt{L_5^3L_2n^{\frac{6\alpha+3d}{8\alpha+2d}}}.
\end{align*} Thus values of $L_2, L_5$ such that \begin{equation*}
    z_1 \sqrt{L_2L_5^3}n^{\frac{6\alpha+3d}{16\alpha+4d}+\frac{3\alpha}{8\alpha+2d}} > 8\max(1,||g_0||_\infty)\log^{9(d+1)/4}(n)
\end{equation*} are sufficient to verify \eqref{QGCP2}. For $n\geq3$, $d>1$ and $\alpha \in [0,1]$ $n^{\frac{12\alpha+3d}{16\alpha+4d}} > \log^{9(d+1)/4}(n)$ so \begin{equation*}
    L_2L_5^3 > \bigg(\frac{8\max(1,||g_0||_\infty)}{(3/L_1)^{3/2}d^{1/4}\sqrt{2\tau}}\bigg)^2
\end{equation*} is a sufficient condition to verify \eqref{QGCP2}.

\subsubsection{Verifying \eqref{QGCP3}}
First consider, \begin{equation*}
    L_1\bar\delta_n \leq 4L_1\max(1,||g_0||_\infty)n^{-\alpha/(4\alpha+d)}\log^{2\rho+2d+2}(n), 
\end{equation*} and \begin{equation*}
    3\sqrt{||\mu||}\beta_n \geq 3\sqrt{||\mu||}(L_5+L_6+L_7)n^{(2\alpha+d)/(8\alpha+2d)}\log^{2\rho+(d+1)/2}(n).
\end{equation*} Plainly $n^{(2\alpha+d)/(8\alpha+2d)}\log^{2\rho+(d+1)/2}(n)>n^{-\alpha/(4\alpha+d)}\log^{2\rho+2d+2}(n)$ for $n \geq 3$, so \begin{equation*}
    L_5+L_6+L_7 > \frac{4L_1\max(1,||g_0||_\infty)}{3\sqrt{||\mu||}}
\end{equation*} is a sufficient condition to verify \eqref{QGCP3}.

\subsubsection{Verifying \eqref{QGCP4}}
Consider, \begin{align*}
    &\quad m\zeta_n^d\bigg(\log\bigg(\frac{3^{3/2}d^{1/4}\beta_n^{3/2}\sqrt{2\tau \zeta_n}}{(L_1\bar\delta_n)^{3/2}} \bigg) \bigg)^{1+d} \\
    &\leq m(L_2^d+L_3^d+L_4^d)n^{\frac{2\alpha+d}{4\alpha+d}}\log^{4\rho}(n)\bigg(\frac{3}{2}\log\Big(\frac{3\beta_n}{L_1\bar\delta_n} \Big)+\frac{1}{2}\log(2\tau d^{1/2}\zeta_n)\bigg)^{1+d}, \\
    &\leq m(L_2^d+L_3^d+L_4^d)n^{\frac{2\alpha+d}{4\alpha+d}}\log^{4\rho}(n)\bigg(\frac{3}{2}\log\bigg(\frac{3(L_5+L_6+L_7)}{2||g_0||_\infty}n^{\frac{6\alpha+d}{8\alpha+2d}}\log^\rho(n)\bigg)+\\
    &\quad \quad \quad \quad \quad \quad \quad \quad \quad \quad \quad \quad \frac{1}{2}\log(2\tau\sqrt{d}(L_2+L_3+L_4)n^{\frac{2\alpha+d}{4\alpha d+d^2}}\log^{4\rho/d}(n)) \bigg)^{1+d} \\
    &\leq \frac{m(L_2^d+L_3^d+L_4^d)}{2^{1+d}}\bigg(3\log\Big(\frac{3(L_5+L_6+L_7)}{2||g_0||_\infty} \Big)+ \frac{18\alpha + 3d}{8\alpha+2d}\log(n)+3\rho\log(\log(n)) \\
    &\quad \quad \quad \quad \quad \quad \quad \quad \quad \quad \quad \quad +\log(2\tau\sqrt{d}(L_2+L_3+L_4))+ \frac{2\alpha + d}{4\alpha d + d^2}\log(n)+ \frac{4\rho}{d}\log(\log(n))\bigg)^{1+d} n^{\frac{2\alpha+d}{4\alpha+d}}\log^{4\rho}(n)\\
    &= \frac{m(L_2^d+L_3^d+L_4^d)}{2^{1+d}}\bigg(\log\Big(\frac{27\tau\sqrt{d}(L_5+L_6+L_7)^3}{4||g_0||_\infty^3}(L_2+L_3+L_4)\Big) + \bigg(\frac{18\alpha d + 4\alpha +2d+3d^2}{8\alpha d+2d^2}\bigg)\log(n) \\
    &\quad \quad \quad \quad \quad \quad \quad \quad \quad \quad \quad \quad + \bigg(3\rho+\frac{4\rho}{d}\bigg)\log(\log(n))\bigg)^{1+d} n^{\frac{2\alpha+d}{4\alpha+d}}\log^{4\rho}(n) \\
    &\leq \frac{m\sum_{i=2}^4 L_i^d}{2^{1+d}}\bigg(\log\Big(\frac{27\tau\sqrt{d}(\sum_{i=5}^7L_i)^3}{4||g_0||_\infty^3}\sum_{i=2}^4 L_i\Big)+ \Big(\frac{18\alpha d + 4\alpha +2d+3d^2}{8\alpha d+2d^2}+\frac{3\rho d+4\rho}{d} \Big)\log(n)\bigg)^{1+d} n^{\frac{2\alpha+d}{4\alpha+d}}\log^{4\rho}(n) \\
    &\leq  \frac{m\sum_{i=2}^4 L_i^d}{2^{1+d}}\bigg(\log\Big(\frac{27\tau\sqrt{d}(\sum_{i=5}^7L_i)^3}{4||g_0||_\infty^3}\sum_{i=2}^4 L_i\Big)+ \frac{32\alpha d + 12\alpha +2d+3d^2+6\alpha d^2}{8\alpha d+2d^2}  \log(n)\bigg)^{1+d} n^{\frac{2\alpha+d}{4\alpha+d}}\log^{4\rho}(n) \\
    &\leq \frac{m\sum_{i=2}^4 L_i^d}{2^{1+d}}\bigg(\log\Big(\frac{27\tau\sqrt{d}(\sum_{i=5}^7L_i)^3}{4||g_0||_\infty^3}\sum_{i=2}^4 L_i\Big)+ \Big(4+\frac{12+d+d^2}{8\alpha d+2d^2}\Big)  \log(n)\bigg)^{1+d} n^{\frac{2\alpha+d}{4\alpha+d}}\log^{4\rho}(n)
\end{align*} and \begin{align*}
    n\bar\delta_n^2 &\geq 4||g_0||_\infty^2 n^{\frac{2\alpha+d}{4\alpha+d}}\log^{2\rho+2(d+1)}(n)
\end{align*} Condition \eqref{QGCP4} is then satisfied for all $n$ such that \begin{align*}
    4||g_0||_\infty^2 \log^{2d+2-2\rho}(n) &\geq  \frac{m\sum_{i=2}^4 L_i^d}{2^{1+d}}\bigg(\log\Big(\frac{27\tau\sqrt{d}(\sum_{i=5}^7L_i)^3 \sum_{i=2}^4 L_i}{4||g_0||_\infty^3}\Big) \\
    &\quad \quad \quad \quad +\Big(4+ \frac{12+d+d^2}{8\alpha d + 2d^2}\Big)\log(n)\bigg)^{1+d}
    \end{align*}

\subsubsection{Verifying \eqref{QGCP5}}
Consider, \begin{align*}
    2\log\bigg(\frac{6\beta_n\sqrt{||\mu||}}{L_1\bar\delta_n}\bigg)&\leq 2\log\bigg(\frac{6\sqrt{||\mu||}(L_5+L_6+L_7)n^{(2\alpha+d)/(8\alpha+2d)}\log^{4\rho+(d+1)/2}(n)}{(2||g_0||_\infty+1) L_1n^{-2\alpha/(4\alpha+d)}\log^{\rho+(d+1)/2}(n)} \bigg) \\
    &= 2\log\bigg(\frac{6\sqrt{||\mu||}(L_5+L_6+L_7)}{2L_1||g_0||_\infty+L_1}n^{(6\alpha+d)/(8\alpha+2d)}\log^{3\rho}(n) \bigg) \\
    &= 2\log\bigg( \frac{6\sqrt{||\mu||}(L_5+L_6+L_7)}{2L_1||g_0||_\infty+L_1}\bigg) + \log\bigg( n^{(6\alpha+d)/(4\alpha+d)}\log^{6\rho}(n) \bigg) 
\end{align*} and \begin{align*}
    n\bar\delta_n^2&=4||g_0||_\infty^2 n^{(2\alpha+d)/(4\alpha+d)}\log^{2\rho+2d+2}(n) + 4||g_0||_\infty n^{(\alpha+d)/(4\alpha+d)}\log^{3\rho+3d+3}(n) + n^{d/(4\alpha+d)}\log^{4\rho+4d+4}(n) \\
    &\geq 4||g_0||_\infty^2 n^{(4\alpha+2d)/(8\alpha+2d)}\log^{2\rho+2d+2}(n) 
\end{align*} 

Therefore, condition \eqref{QGCP5} holds for all $n$ such that \begin{equation*}
    2\log\bigg( \frac{6\sqrt{||\mu||}(L_5+L_6+L_7)}{2L_1||g_0||_\infty+L_1}\bigg) \leq 4||g_0||_\infty^2 n^{(4\alpha+2d)/(8\alpha+2d)}\log^{2\rho+2d+2}(n) - \log\bigg( n^{(6\alpha+d)/(4\alpha+d)}\log^{6\rho}(n) \bigg) 
\end{equation*}

\subsubsection{Verifying \eqref{QGCP6}}

Consider \begin{align*}
     \beta_n^2 &= L_5^2 n^{\frac{2\alpha+d}{4\alpha+d}}\log^{4\rho+d+1}(n) + L_6^2n^{\frac{\alpha+d}{4\alpha+d}}\log^{6\rho+d+1}(n) + L_7^2 n^{\frac{d}{4\alpha+d}}\log^{8\rho+d+1}(n) 
\end{align*} and \begin{align*}
    \zeta_n^d \bigg(\log\bigg(\frac{3\zeta_n}{\bar\delta_n} \bigg) \bigg)^{1+d} &= ( L_2^d n^{\frac{2\alpha+d}{4\alpha+d}}\log^{2\rho}(n) + L_3^dn^{\frac{\alpha+d}{4\alpha+d}}\log^{3\rho}(n) + L_4^d n^{\frac{d}{4\alpha+d}}\log^{4\rho}(n)) \\
    &\quad \times \bigg(\log \bigg(\frac{3L_2n^{(2\alpha+d)/(4\alpha d+d^2)}\log^{2\rho/d}(n) + 3L_3n^{(\alpha+d)/(4\alpha d+d^2)}\log^{3\rho/d}(n) +3L_4n^{d/(4\alpha d+d^2)}\log^{4\rho/d}(n)}{2||g_0||_\infty n^{-\alpha/(4\alpha+d)}\log^{\rho+d+1}(n) + n^{-2\alpha/(4\alpha+d)}\log^{2\rho+2d+2}(n)} \bigg) \bigg)^{1+d} \\
    &\leq ( L_2^d n^{\frac{2\alpha+d}{4\alpha+d}}\log^{2\rho}(n) + L_3^dn^{\frac{\alpha+d}{4\alpha+d}}\log^{3\rho}(n) + L_4^d n^{\frac{d}{4\alpha+d}}\log^{4\rho}(n)) \\
    &\quad \times \bigg( \log \bigg(\frac{3(L_2+L_3+L_4)n^{(2\alpha+d)/(4\alpha d+d^2)}\log^{4\rho/d}(n)}{\min(1,2||g_0||_\infty) n^{-2\alpha/(4\alpha+d)}\log^{\rho+d+1)}(n)} \bigg)\bigg)^{1+d} \\
    &\leq  ( L_2^d n^{\frac{2\alpha+d}{4\alpha+d}}\log^{2\rho}(n) + L_3^dn^{\frac{\alpha+d}{4\alpha+d}}\log^{3\rho}(n) + L_4^d n^{\frac{d}{4\alpha+d}}\log^{4\rho}(n)) \\
    &\quad \times \Bigg[\log\bigg(\frac{3(L_2+L_3+L_4)}{\min(1,2||g_0||_\infty)} \bigg) + \log\Big(n^{\frac{2\alpha+d}{4\alpha d + d^2}+\frac{2\alpha}{4\alpha+d}}\log^{4\rho-\rho/d-d-1}(n) \Big)\Bigg]^{1+d}.
\end{align*} Define \begin{align*}
    \mathcal{K}(n) &= \Bigg[\log\bigg(\frac{3(L_2+L_3+L_4)}{\min(1,2||g_0||_\infty)} \bigg) + \frac{2\alpha d+ 2\alpha +d}{4\alpha d+d^2}\log(n)+(4\rho-\rho/d-d-1)\log(\log(n))\Bigg]^{1+d}\\
    &\leq \log^{1+d}(n)\bigg(\log\bigg(\frac{3(L_2+L_3+L_4)}{\min(1,2||g_0||_\infty)} \bigg) + \frac{2\alpha d+ 2\alpha +d}{4\alpha d+d^2} + (4\rho-\rho/d-d-1) \bigg)^{1+d},
\end{align*} for $n\geq 3$. Let $\mathcal{K}_1=\bigg(\log\bigg(\frac{3(L_2+L_3+L_4)}{\min(1,2||g_0||_\infty)} \bigg) + \frac{2\alpha d+ 2\alpha +d}{4\alpha d+d^2} + (4\rho-\rho/d-d-1) \bigg)$. Grouping terms of the same order we require the following for all sufficiently large $n$ \begin{align*}
    L_5^2 \log^{2\rho}(n) &\geq 16K_5L_2^d \mathcal{K}_1^{1+d}, \\
    L_6^2 \log^{3\rho}(n) &\geq 16K_5L_3^d \mathcal{K}_1^{1+d}, \\
    L_7^2 \log^{4\rho}(n) &\geq 16K_5L_4^3 \mathcal{K}_1^{1+d},
\end{align*} to satisfy \eqref{QGCP6}. Thus, the following are sufficient conditions to satisfy \eqref{QGCP6} for all $n \geq 3$ \begin{align*}
    L_5 \geq \sqrt{\frac{16K_5L_2^d\mathcal{K}_1^{1+d}}{\log^{2\rho}(3)}}, \quad     L_6 \geq \sqrt{\frac{16K_5L_3^d\mathcal{K}_1^{1+d}}{\log^{3\rho}(3)}}, \quad     L_7 \geq \sqrt{\frac{16K_5L_4^d\mathcal{K}_1^{1+d}}{\log^{4\rho}(3)}}.
\end{align*}

\subsubsection{Verifying \eqref{QGCP7}}
Consider, \begin{align*}
    2c_5n\delta_n^2 &=2c_5\bigg(4||g_0||_\infty^2 n^{\frac{2\alpha+d}{4\alpha+d}}\log^{2\rho}(n) + 4||g_0||_\infty n^{\frac{\alpha+d}{4\alpha+d}}\log^{3\rho}(n)+ n^{\frac{d}{4\alpha+d}}\log^{4\rho}(n)\bigg),
\end{align*} and \begin{align*}
    D_1\zeta_n^d \log^{q_2}(\zeta_n) \geq D_1\bigg(L_2 n^{\frac{2\alpha+d}{4\alpha+d}}\log^{2\rho}(n) + L_3 n^{\frac{\alpha+d}{4\alpha+d}}\log^{3\rho}(n)+ L_4n^{\frac{d}{4\alpha+d}}\log^{4\rho}(n)\bigg)
\end{align*} for $\zeta_n>e$ - i.e. such that $\log(\zeta_n) \geq 1$. Then condition \eqref{QGCP1} and $L_2 \geq (8c_5||g_0||_\infty^2)/D_1$, $L_3 \geq (8c_5||g_0||_\infty)/D_1$ and $L_4 \geq 2c_5/D_1$ are sufficient conditions to verify \eqref{QGCP7}.
\subsubsection{Verifying \eqref{QGCP8}}
Consider, \begin{align*}
    8c_5n\delta_n^2 &=8c_5\bigg(4||g_0||_\infty^2 n^{\frac{2\alpha+d}{4\alpha+d}}\log^{2\rho}(n) + 4||g_0||_\infty n^{\frac{\alpha+d}{4\alpha+d}}\log^{3\rho}(n)+ n^{\frac{d}{4\alpha+d}}\log^{4\rho}(n)\bigg),
\end{align*} and \begin{align*}
    \beta_n^2 \geq L_5^2 n^{\frac{2\alpha+d}{4\alpha+d}}\log^{2\rho}(n) + L_6^2n^{\frac{\alpha+d}{4\alpha+d}}\log^{3\rho}(n) + L_7^2 n^{\frac{d}{4\alpha+d}}\log^{4\rho}(n).
\end{align*} Then \eqref{QGCP7} is satisfied with $L_5^2>32||g_0||_\infty^2c_5$, $L_6^2>32||g_0||_\infty c_5$, and $L_7^2>8c_5$.

\subsubsection{Verifying \eqref{QGCP9}}
Consider \begin{align*}
    \exp(c_5n\delta_n^2)&= \exp\Bigg(c_5\bigg(4||g_0||_\infty^2 n^{\frac{2\alpha+d}{4\alpha+d}}\log^{2\rho}(n) + 4||g_0||_\infty n^{\frac{\alpha+d}{4\alpha+d}}\log^{3\rho}(n)+ n^{\frac{d}{4\alpha+d}}\log^{4\rho}(n)\bigg)\Bigg), \\
    &\geq \exp\bigg(4c_5||g_0||^2_\infty n^{(2\alpha+d)/(4\alpha +d)}\log^{2\rho}(n)\bigg)
\end{align*} and \begin{align*}
    \zeta_n^{q_1-d+1}&\leq \bigg(L_2n^{(2\alpha+d)/(4\alpha d+d^2)}\log^{2\rho/d}(n) + L_3n^{(\alpha+d)/(4\alpha d+d^2)}\log^{3\rho/d}(n) +L_4n^{d/(4\alpha d+d^2)}\log^{4\rho/d}(n)\bigg)^{q_1}, \\
    &\leq \bigg((L_2+L_3+L_4)n^{(2\alpha+d)/(4\alpha d + d^2)}\log^{4\rho/d}(n)\bigg)^{q_1}
\end{align*} The condition is then satisfied for all $n$ such that \begin{align*}
    n^{(2\alpha+d)/(4\alpha +d)}\log^{2\rho}(n) \geq \frac{q_1}{4c_5||g_0||^2_\infty}\log\Big((L_2+L_3+L_4)n^{(2\alpha+d)/(4\alpha d + d^2)}\log^{4\rho/d}(n) \Big).
\end{align*}

\subsection{SGCP model}
Recall the definitions of the following sequences: \begin{align*}
    \delta_n &= n^{-\alpha/(2\alpha +d}\log^{\rho}(n)\\
    \bar\delta_n &= n^{-\alpha/(2\alpha +d}\log^{\rho+d+1}(n) \\
    \zeta_n &= L_8n^{\frac{1}{2\alpha+d}}(\log(n))^{2\rho/d}, \\
    \beta_n &= L_9n^{\frac{d}{2(2\alpha+d)}}(\log(n))^{d+1+2\rho}, \\
    \lambda_n &= L_{10}n^{\frac{d}{\kappa(2\alpha+d)}}(\log(n))^{4\rho/\kappa}
\end{align*}  with $L_8, L_9, L_{10}$ satisfying \begin{align*}
L_8 &> \max\Bigg( A, 1, \bigg(\frac{2c_5}{D_1} \bigg)^{1/d} \Bigg) \\
L_9 &\geq \sqrt{8c_5} \\
L_{10} &> \bigg( \frac{c_5}{c_0}\bigg)^{1/\rho} \\
L_8L_9^3L_{10}^{3/2} &> \frac{2}{()6cL_1)^{3/2}\tau\sqrt{d}} \\
L_9L_{10}^{1/2} &> \frac{1}{6cL_1\sqrt{||\mu||}}
\end{align*}for $n\geq \max(3,n_6,n_7,n_8)$. Here $n_6$ is the smallest integer such that \begin{equation*}
n^{\frac{d}{2\alpha+d}} > \max \Bigg( 2\log(12cL_1L_9L_{10}^{1/2})+1, \log(2L_1L_{10}^{1/2})+1, \frac{1}{c_5}\big(\log(L_8^{q_1-d+1})+1 \big) \Bigg),
\end{equation*} $n_7$ is the smallest integer such that \begin{equation*}
\log^{2d+2}(n) >    mL_8^d \bigg(\log((6cL_1)^{3/2}\sqrt{2\tau L_8 L_9^3}L_{10}^{3/4}d^{1/4})+\frac{\kappa(6d+6\alpha+2)+3d}{4\kappa(2\alpha+d)}\log(n)+\log\Big(\log^{3\rho/2 +3\rho/\kappa+\rho/d-d-1}(n)\Big)\bigg)^{1+d},
\end{equation*} and $n_8$ is  the smallest integer such that \begin{equation*}
\log^{2\rho}(n) > \frac{16K_5D_1L_8^d}{L_9^2}\Bigg(\frac{\log(\sqrt{L_{10}}L_8)+\log\big(n^{\frac{2\alpha\kappa+ 2\kappa+d}{2\kappa(2\alpha+d)}}\log^{\rho(2/\kappa+2/d-1)-d-1}(n)\big) }{\log(n)} \Bigg)^{1+d}.
\end{equation*}

In the remainder of this subsection, we show that the following conditions, which are all restatements of required results in our main analysis, hold for the sequences described above. \begin{align}
    (6cL_1)^{3/2}d^{1/4}\beta_n^{3/2}\lambda_n^{3/4}\sqrt{2\tau\zeta_n} &> 2\bar\delta_n^{3/2} \label{SGCP1} \\
   6cL_1 \beta_n\sqrt{\lambda_n||\mu||} &> \bar\delta_n \label{SGCP2} \\
    \zeta_n &> \max(A,1) \label{SGCP3} \\
    m \zeta_n^d\bigg(\log\Big(\frac{(6cL_1)^{3/2}\lambda_n^{3/4}\beta_n^{3/2}d^{1/4}\sqrt{2\tau\zeta_n}}{\bar\delta_n^{3/2}} \Big) \bigg)^{1+d} &< K_3n\bar\delta_n^2 \label{SGCP4} \\
    2 \log\Big(\frac{12cL_1\beta_n\sqrt{\lambda_n||\mu||}}{\bar\delta_n} \Big) &< K_4n\bar\delta_n^2 \label{SGCP5} \\
    \log\Big( \frac{2L_1\lambda_n^{1/2}}{\bar\delta_n}\Big) &< K_5 n \bar\delta_n^2 \label{SGCP6} \\
    \beta_n^2 &> 16K_5\zeta_n^d\bigg(\log\Big(\frac{\lambda_n^{1/2}\zeta_n}{\bar\delta_n}\Big) \bigg)^{1+d} \label{SGCP7} \\
    c_0\lambda_n^\rho &>c_5n\delta_n^2 \label{SGCP8} \\
    D_1\zeta_n^d &\geq 2c_5n\delta_n^2 \label{SGCP9} \\
    \zeta^{q_1-d+1}_n &\leq e^{c_5n\delta_n^2} \label{SGCP10} \\
    \beta_n^2 &\geq 8c_5n\delta_n^2 \label{SGCP11}
\end{align} 

In turn we demonstrate that each of the conditions \eqref{SGCP1} through \eqref{SGCP11} hold.

\subsubsection{Verifying \eqref{SGCP1}}
Consider \begin{align*}
    &\enspace (6cL_1)^{3/2}d^{1/4}\beta_n^{3/2}\lambda_n^{3/4}\sqrt{2\tau\zeta_n} \\
    &= (6cL_1)^{3/2}d^{1/4}L_9^{3/2}n^{\frac{3d}{4(2\alpha+d)}}\log^{\frac{3d+3}{4}+3\rho}(n)L_{10}^{3/4}n^{\frac{3d}{4\kappa(2\alpha+d)}}\log^{3\rho/\kappa}(n)\sqrt{2\tau L_8n^{\frac{1}{2\alpha+d}}(\log(n))^{2\rho/d}} \\
    &=(6cL_1)^{3/2}\sqrt{2\tau L_8 L_9^3}L_{10}^{3/4}d^{1/4} n^{\frac{3d}{4(2\alpha+d)}+\frac{3d}{4\kappa(2\alpha+d)}+\frac{1}{2(2\alpha+d)}}\log^{\frac{3d+3}{4}+3\rho +3\rho/\kappa+\rho/d}(n)
\end{align*} and \begin{align*}
    2\bar\delta_n^{3/2}&=2n^{\frac{-3\alpha}{2(2\alpha+d)}}\log^{3\rho/2+3(d+1)/2}(n)
\end{align*}
So \eqref{SGCP1} can be rewritten: \begin{equation*}
    (6cL_1)^{3/2}\sqrt{2\tau L_8 L_9^3}L_{10}^{3/4}d^{1/4} n^{\frac{3d}{4(2\alpha+d)}+\frac{3d}{4\kappa(2\alpha+d)}+\frac{3\alpha+1}{2(2\alpha+d)}}\log^{3\rho/2 +3\rho/\kappa+\rho/d-\frac{3d-3}{4}}(n) > 2,
\end{equation*} which holds for $L_8, L_9, L_{10}$ such that $L_8 L_9^3 L_{10}^{3/2} > \frac{2}{(6cL_1)^{3/2}\tau\sqrt{d}}$.

\subsubsection{Verifying \eqref{SGCP2}}
We may rewrite \eqref{SGCP2} as \begin{equation*}
   6cL_1 \sqrt{||\mu||}L_9L_{10}^{1/2}n^{\frac{1}{2}}\log^{\rho+2\rho/\kappa-d-1}(n) >1
\end{equation*} which holds for all $L_9,L_{10}$ such that $L_9L_{10}^{1/2}> 1/(6cL_1\sqrt{||\mu||})$.

\subsubsection{Verifying \eqref{SGCP3}}
If $n\geq 3$ then $\zeta_n$ holds for all $L_8 \geq \max(A,1)$.

\subsubsection{Verifying \eqref{SGCP4}}
Consider \begin{align*}
    &m \zeta_n^d\bigg(\log\Big(\frac{(6cL_1)^{3/2}\lambda_n^{3/4}\beta_n^{3/2}d^{1/4}\sqrt{2\tau\zeta_n}}{\bar\delta_n^{3/2}} \Big) \bigg)^{1+d} \\
    &= m L_8^dn^{\frac{d}{2\alpha +d}}\log^{2\rho}(n) \log\bigg((6cL_1)^{3/2}\sqrt{2\tau L_8 L_9^3}L_{10}^{3/4}d^{1/4} n^{\frac{3d}{4(2\alpha+d)}+\frac{3d}{4\kappa(2\alpha+d)}+\frac{3\alpha+1}{2(2\alpha+d)}}\log^{3\rho/2 +3\rho/\kappa+\rho/d-d-1}(n)\bigg)^{1+d}
\end{align*} and \begin{align*}
    n\bar\delta_n^2 &= n^{\frac{d}{2\alpha+d}}\log^{2\rho+2d+2}(n) 
\end{align*} Thus \eqref{SGCP4} holds for all $n$ such that \begin{align*}
  \log^{2d+2}(n) >    mL_8^d \bigg(\log((6cL_1)^{3/2}\sqrt{2\tau L_8 L_9^3}L_{10}^{3/4}d^{1/4})+\frac{\kappa(6d+6\alpha+2)+3d}{4\kappa(2\alpha+d)}\log(n)+\log\Big(\log^{3\rho/2 +3\rho/\kappa+\rho/d-d-1}(n)\Big)\bigg)^{1+d}
\end{align*}

\subsubsection{Verifying \eqref{SGCP5}}
We may rewrite \eqref{SGCP5} as \begin{equation*}
    2\log\bigg( 12cL_1L_9 L_{10}^{1/2}n^{\frac{1}{2}+\frac{d}{2\kappa(2\alpha+d)}}\log^{2\rho/\kappa+\rho-d-1}(n)\bigg) <  n^{\frac{d}{2\alpha+d}}\log^{2\rho+d+1}(n) 
\end{equation*} which holds for all $n$ such that \begin{equation*}
n^{\frac{d}{2\alpha+d}}> 2\log(12cL_1L_9L_{10}^{1/2})+1.
\end{equation*}

\subsubsection{Verifying \eqref{SGCP6}}
We may rewrite \eqref{SGCP6} as \begin{equation*}
    \log\bigg(2L_1L_{10}^{1/2}n^{\frac{d}{2\kappa(2\alpha +d)}+\frac{\alpha}{2\alpha+d}}\log^{2\rho/\kappa-\rho-d-1}(n) \bigg) < n^{\frac{d}{2\alpha+d}}\log^{2\rho+d+1}(n)
\end{equation*} which holds for all $n$ such that \begin{equation*}
n^{\frac{d}{2\alpha+d}}> \log(2L_1L_{10}^{1/2})+1.
\end{equation*}

\subsubsection{Verifying \eqref{SGCP7}}
Consider \begin{align*}
    \beta_n^2 = L_9^2n^{\frac{d}{2\alpha +d}}\log^{d+1+4\rho}(n)
\end{align*} and \begin{align*}
    16K_5\zeta_n^d \bigg(\log\Big(\frac{\lambda_n^{1/2}\zeta_n}{\bar\delta_n}\Big) \bigg)^{1+d} &= 16K_5D_1L_8^dn^{\frac{d}{2\alpha +d}}\log^{2\rho}(n)\bigg(\log\Big( \frac{\sqrt{L_{10}}L_8 n^{\frac{2\kappa+d}{2\kappa(2\alpha+d)}}\log^{2\rho/\kappa2 + 2\rho/d}(n)}{n^{\frac{-\alpha}{2\alpha+d}}\log^{\rho+d+1}(n)}\Big) \bigg)^{1+d} \\
    &= 16K_5D_1L_8^dn^{\frac{d}{2\alpha +d}}\log^{2\rho}(n)\bigg(\log(\sqrt{L_{10}}L_8)+\log\Big(n^{\frac{2\alpha\kappa+ 2\kappa+d}{2\kappa(2\alpha+d)}}\log^{\rho(2/\kappa+2/d-1)-d-1}(n) \bigg)^{1+d}
\end{align*}

So \eqref{SGCP7} can then be rewritten as \begin{equation*}
    L_9^2\log^{d+1+2\rho}(n) > 16K_5 D_1 L_8^d \bigg(\log(\sqrt{L_{10}}L_8)+\log\Big(n^{\frac{2\alpha\kappa+ 2\kappa+d}{2\kappa(2\alpha+d)}}\log^{\rho(2/\kappa+2/d-1)-d-1}(n)\Big) \bigg)^{1+d}
\end{equation*} which holds for all $n$ such that \begin{equation*}
\log^{2\rho}(n) > \frac{16K_5D_1L_8^d}{L_9^2}\Bigg(\frac{\log(\sqrt{L_{10}}L_8)+\log\big(n^{\frac{2\alpha\kappa+ 2\kappa+d}{2\kappa(2\alpha+d)}}\log^{\rho(2/\kappa+2/d-1)-d-1}(n)\big) }{\log(n)} \Bigg)^{1+d}
\end{equation*}

\subsubsection{Verifying \eqref{SGCP8}}
We may rewrite \eqref{SGCP8} as \begin{align*}
    c_0L_{10}^\rho n^{\frac{d\rho}{\kappa(2\alpha+d)}}\log^{4\rho^2/\kappa}(n) > c_5n^{\frac{d}{2\alpha+d}}\log^{2\rho}(n)
\end{align*} If $\rho/\kappa>1$ this holds for all $L_{10}>(c_5/c_0)^{1/\rho}$.

\subsubsection{Verifying \eqref{SGCP9}}
We may rewrite \eqref{SGCP9} as \begin{align*}
    D_1L_8^dn^{\frac{d}{2\alpha +d}}\log^{2\rho}(n) > 2c_5n^{\frac{d}{2\alpha+d}}\log^{\frac{2+2d}{2+d/\alpha}}(n)
\end{align*} which is satisfied for all $L_8 > (2c_5/D_1)^{1/d}$.

\subsubsection{Verifying \eqref{SGCP10}}
Consider \begin{align*}
    \zeta_n^{q_1-d+1}= L_8^{q_1-d+1} n^{\frac{q_1-d+1}{2\alpha+d}}\log^{\frac{2\rho(q_1-d+1)}{d}}(n)
\end{align*} and \begin{align*}
    \exp(c_5n\delta_n^2)= \exp\big(c_5 n^{\frac{d}{2\alpha+d}}\log^{2\rho}(n) \big).
\end{align*} Then \eqref{SGCP10} holds for all $n$ such that \begin{equation*}
n^{\frac{d}{2\alpha+d}}> \frac{1}{c_5}\big(\log(L_8^{q_1-d+1})+1 \big).
\end{equation*}

\subsubsection{Verifying \eqref{SGCP11}}
We may rewrite \eqref{SGCP11} as \begin{align*}
    L_9^2n^{\frac{d}{2\alpha +d}}\log^{d+1+4\rho}(n) \geq 8c_5n^{\frac{d}{2\alpha+d}}\log^{2\rho}(n).
\end{align*} which is satisfied for all $L_9\geq\sqrt{8c_5}$.

\section*{Acknowledgements}
The authors thank Marco Battison for helpful discussions during the preparation of this paper.

\bibliographystyle{apalike} 
\bibliography{GCPrefs}

\end{document}